\documentclass[11pt,a4paper]{article} 

   \usepackage{amsmath,amscd,amssymb,latexsym,graphicx,color} 
 
\newtheorem{defn}{Definition} 
\newtheorem{lem}{Lemma} 
\newtheorem{prop}{Proposition} 
\newtheorem{thm}{Theorem}

\newtheorem{rem}{Remark} 
\newtheorem{rems}{Remarks}

\numberwithin{equation}{section}

\reversemarginpar 
 
\newcommand{\op}{\text{op}} 
\newcommand{\iso}{\text{Iso}}

\input utile.def 

\newcommand{\ot}{\otimes}
\newcommand{\hot}{\hat{\otimes}} 

\newcommand{\To}{\rightarrow}

\newcommand{\Tb}{\overline{T^qX}} 
\newcommand{\ev}[1][]{\mbox{ev}_{#1}} 
\newcommand{\Ind}[2][]{\mbox{Ind}^{#2}(#1)}

\newcommand{\Id}{\mbox{Id}} 

\newcommand{\Exp}{\mbox{exp}}
\newcommand{\dom}{\mbox{dom}}
\newcommand{\interior}{\mbox{int}}

\newcommand{\ba}{\mathbf{a}}
\newcommand{\fraq}[2]{\displaystyle{\frac{#1}{#2}}}

\title{Elliptic symbols, elliptic operators and Poincar\'e duality on conical pseudomanifolds}  
 \author{Jean-Marie Lescure}
 
\begin{document} 
 
\maketitle

{\small Abstract: In \cite{DL}, a notion of noncommutative tangent
  space is associated with a conical pseudomanifold and the Poincar\'e duality in $K$-theory
  is proved between this space and the pseudomanifold. The present
  paper continues this work. We show that an appropriate presentation of
  the notion of symbols on a manifold generalizes right away to
  conical pseudomanifolds and that it enables us to interpret the
  Poincar\'e duality in the singular setting as a noncommutative symbol map.}

\tableofcontents

\section{Introduction}

In this paper we give a concrete description of the Poincar\'e duality
in $K$-theory for a conical pseudomanifold as stated and
proved in \cite{DL}. This duality holds between the algebra $C(X)$ of continuous fonctions on a
(compact) pseudomanifold $X$  and the $C^*$-algebra $C^*(T^cX)$ of a suitable
 {\sl tangent space} of this pseudomanifold.  

\smallskip The tangent space $T^cX$ introduced in \cite{DL} is a smooth groupoid. It 
is no more commutative, but it restricts  to the usual
tangent space of a manifold outside the singularity and the singular
contribution is quite simple.

\smallskip The duality between $C(X)$ and $C^*(T^cX)$ is defined in terms of bivariant $K$-theory
but it is important to recall that it implies the existence
of an isomorphism: 
 \begin{equation}\label{PD-map-conic-case} 
  \Sigma^c\ : \  K_0(X)\overset{\simeq}{\longrightarrow} K_0(C^*(T^cX))
 \end{equation}

\smallskip The main purpose of this paper is to identify this isomorphism with a
noncommutative symbol map, as one does in the smooth case with the
usual symbol map. Indeed, the Poincar\'e duality in the case of a smooth
closed manifold $V$ induces an isomorphism between $K_0(V)$ and
$K_0(C^*(TV))\simeq K^0(T^*V)$ which is nothing else but the principal
symbol map: 
 \begin{equation}\label{PD-map-smooth-case} 
    \begin{matrix} 
          K_0(V) & \longrightarrow & K^0(T^*V) \\
             [P] & \longmapsto    & [\sigma(P)] 
    \end{matrix}
 \end{equation}
sending classes of elliptic pseudodifferential operators (the basic cycles of the
$K$-homology of $V$) to classes of their principal symbols (the basic
cycles of the $K$-theory with compact supports of $T^*V$).  

\smallskip The interpretation of (\ref{PD-map-conic-case}) as a
noncommutative symbol map is
really important for two reasons. Firstly, it validates the
choice of a tangent space which is $K$-dual to the singular manifold
and thus motivates further investigations toward analysis or
differential geometry of singular spaces by using this noncommutative
object as well as other tools of noncommutative geometry. 
Secondly, this approach can be of interest
for people looking for  Fredholmness conditions in elliptic equations in
singular situations like stratified spaces. Indeed, the
notion of tangent space of a stratified space is very intuitive
as soon as one understands the conical case, and the notion of
elliptic noncommutative symbols
appears directly. The case of general stratifications will be treated
in forthcoming articles. 

\smallskip In \cite{DL}, we propose two $KK$--equivalent definitions of the tangent space of a
pseudomanifold $X$ and the main results were stated for the first one,
noted $T^cX$ in the present article. To explain in what sense cycles of the $K$-theory
of the tangent space of a pseudomanifold are {\sl noncommutative symbols} and cycles of its $K$-homology
are {\sl pseudodifferential operators}, we will use here the
second definition given in \cite{DL}, noted in the sequel $T^qX$. The equivalence in $K$-theory of both
tangent spaces allows us to state all the results of \cite{DL} for $T^qX$ and in
particular the isomorphism (\ref{PD-map-conic-case}). Even if this
equivalence is obvious for people familiar with groupoids, one
will give full details about it in section \ref{section2}.

\smallskip Now, surprisingly, one can define {\sl noncommutative symbols}  on a pseudomanifold exactly as
one defines symbols on a smooth manifold. More precisely, symbols on a smooth manifold $V$ are functions on the
cotangent space $T^*V$ with adequate behavior in the fibers.  They can
be considered as pointwise multiplication operators on, for instance,
$C^\infty_c(T^*V)$. Under a Fourier transform  in the fibers, they can
also be viewed as families parametrized by $V$ of convolutions
operators in the fibers of $TV$. Thus: 

\centerline{Symbols on $V$ are
pseudodifferential operators on the tangent space $TV$,}

\noindent  where $TV$ is considered as a {\sl} groupoid and we talk
about pseudodifferential calculus for groupoids \cite{NWX,MP,Vas,Vas2}. 

\smallskip This simple observation is already important to understand that
the tangent groupoid defined by A. Connes in \cite{co0} gives the
analytic index of elliptic pseudodifferential operators. Next, it
suggests the following definition of noncommutative symbols on the pseudomanifold
$X$. 

\centerline{Noncommutative symbols on $X$ are  pseudodifferential operators on the
tangent space $T^qX$.} 

\smallskip We will see that, after some technical precautions on the Schwartz
kernels and on the behavior near the ``end'' of $T^qX$ of these
pseudodifferential operators, this apparently naive idea works. 
For instance, one can recover in a single
object the notions of interior and conormal symbols arising in
boundary problems and the notion of full ellipticity is quite
immediate here. 

\smallskip Concerning the operators involved in the description of the
Poincar\'e duality, some freedom is allowed: basically, all
calculi based on the work of R. Melrose \cite{Mel} ($b$ or $c$ calculi for
instance) as well as on the work of B.W. Schulze
\cite{schulze1,schulze2} can be used
indifferently and lead to various representants of the same
$K$-homology class (that is, to the Poincar\'e dual of a given
elliptic noncommutative symbol). 

\smallskip The main tools used in this paper are Lie groupoids
(see \cite{DL} and the corresponding bibliography), pseudodifferential
calculus (see \cite{Shub,Mel,MP,NWX,Vas2}) bivariant
$K$-theory (see \cite{Ka1,Ka2,Ska1,CoS,Bla,WO,BJ}). 

\smallskip The author mentions that different techniques have been precedently used to
produce results close from the present work by A. Savin (\cite{Sav1},
see also joint works by V. Nazaikinskii, A. Savin and B. Sternin \cite{NSS1,NSS2}).



\subsection{Reviews and Notations}\label{notations}
The range and source maps of groupoids are noted $r$ and $s$. If $A$
is a subset of the space of units $G^{(0)}$ of a groupoid $G$ then
$G|_A$ denotes the subgroupoid $G|_A=r^{-1}(A)\cap s^{-1}(A)$.
All groupoids in the sequel  are smooth (Lie groupoids), endowed with Haar systems in
order to define their $C^*$-algebras. Moreover, they are
amenable (as continuous fields of amenable groupoids \cite{ARe}). 
In particular, there is no ambiguity  about their
$C^*$-algebras and notations for their $K$-theory will be
shortened:
$$ 
  K^i(G):= K_i(C^*(G)) \text{ and }
  KK(G_1,G_2):= KK(C^*(G_1),C^*(G_2)) 
$$
If $f$ is a homomorphism between two $C^{*}$-algebras $A, B$,  the
corresponding class in $KK(A, B)$ will be denoted by $[f]$. 

When a vector bundle $E\to G^{(0)}$ is given, we define a
$C^*(G)$-Hilbert module noted $C^*(G,E)$ by taking the completion of 
$C^\infty_c(G,r^*E)$ for the norm associated with 
the $C^*(G)$-valued product : 
  $$
    <f , g> (\gamma) = \int_{\eta\in G^{r(\gamma)}}
    <f(\eta^{-1}),g(\eta^{-1}\gamma)>_{s(\eta)}.
  $$
We shall use various deformation groupoids
$G=G_1\times\{t=0\}\cup G_2\times]0,1]_{t}$. The restriction morphism
$\ev[t=0]: G\to G_1$ at $t=0$ gives an exact sequence:
\begin{equation}\label{exact-seq-def-gpd}
 0\To C^*(G_2\times]0,1]) \To C^*(G)\overset{\ev[t=0]}{\To} C^*(G_1)\To
 0
\end{equation}
whose ideal is contractible in $KK$-theory. If $G_1$ is
amenable (which will always be the case in this paper), one gets that 
$[\ev[t=0]]\in KK(G,G_2)$ is invertible. The {\sl deformation element}
associated with the deformation groupoid is the Kasparov element defined by 
\begin{equation}\label{def-elt}
  \partial_G = [\ev[t=0]]^{-1}\otimes[\ev[t=1]]\in KK(G_1,G_2)
\end{equation} 
For convenience, the pair groupoid on a set $E$ will be denoted by
$\cC_E$ . 

\smallskip The (open) cone over a space $L$ is  the quotient space 
$cL = (L\times [0,+\infty[)/ L\times\{0\}$.A conical pseudomanifold is a compact metrisable
space $X$ equipped with the following data. There is one singular
point (but everything in the sequel can be written for a finite number) 
which means that a point $c\in X$ is given and that
$X^{o}:=X\setminus\{c\}$ is a manifold. Moreover there is
   an open neighborhood $\cN$ of $c$, a smooth manifold $L$, continuous
maps $h : \cN\to [0,+\infty[$ and $\varphi_c:\cN\to cL$ satisfying the following:
\begin{itemize} 
 \item  $h$ is surjective, $h^{-1}\{0\}=\{c\}$ and  $h : \cN\setminus\{c\}\to
   ]0,+\infty[$ is a smooth submersion, 
 \item  $\varphi_c:\cN\to cL$ is a homeomorphism, smooth outside $c$, such that: 
  $$ 
    \begin{CD} 
      \cN @>\varphi_c>> cL \\
      @VhVV    @V\overline{p_2}VV \\
      [0,+\infty[ @>=>>[0,+\infty[
    \end{CD}
  $$
commutes. Here $\overline{p_2}$ denotes the quotient map of the
second projection $L\times\RR_+\to \RR_+$. 
\end{itemize}

\smallskip Conical pseudomanifolds are the simplest examples of a stratified
space \cite{BHS}. We distinguish two parts in the
regular stratum $X^{o}$: 
 $$ X^{o}= X_-\cup X_+$$ 
where $X_-=h^{-1}]0,1[$, and $X_+=X\setminus X_-$ is a smooth
compact manifold with boundary, the latter being identified with $L$. The identification 
$X_-\simeq]0,1[\times L$ given by $\varphi_c$ will be often used
without mention. The compactification $M=\overline{X^{o}}$ of
$X^{o}$ into a manifold with boundary $L$ will be sometimes
useful. The following picture illustrates the notations just defined:

\vspace{1cm}

\begin{picture}(0,0)%
\includegraphics{cone2.pstex}%
\end{picture}%
\setlength{\unitlength}{1737sp}%
\begingroup\makeatletter\ifx\SetFigFont\undefined%
\gdef\SetFigFont#1#2#3#4#5{%
  \reset@font\fontsize{#1}{#2pt}%
  \fontfamily{#3}\fontseries{#4}\fontshape{#5}%
  \selectfont}%
\fi\endgroup%
\begin{picture}(17343,4383)(586,-3595)
\put(3601,-3361){\makebox(0,0)[lb]{\smash{{\SetFigFont{8}{9.6}{\rmdefault}{\mddefault}{\updefault}{\color[rgb]{0,0,0}$h$}%
}}}}
\put(601,-961){\makebox(0,0)[lb]{\smash{{\SetFigFont{8}{9.6}{\rmdefault}{\mddefault}{\updefault}{\color[rgb]{0,0,0}$c$}%
}}}}
\put(1501,-661){\makebox(0,0)[lb]{\smash{{\SetFigFont{8}{9.6}{\rmdefault}{\mddefault}{\updefault}{\color[rgb]{0,0,0}$L$}%
}}}}
\put(3001,-1711){\makebox(0,0)[lb]{\smash{{\SetFigFont{12}{14.4}{\rmdefault}{\mddefault}{\updefault}{\color[rgb]{0,0,0}$X$}%
}}}}
\put(5701,-811){\makebox(0,0)[lb]{\smash{{\SetFigFont{8}{9.6}{\rmdefault}{\mddefault}{\updefault}{\color[rgb]{0,0,0}$L$}%
}}}}
\put(751,-3361){\makebox(0,0)[lb]{\smash{{\SetFigFont{8}{9.6}{\rmdefault}{\mddefault}{\updefault}{\color[rgb]{0,0,0}$0$}%
}}}}
\put(8551,-3361){\makebox(0,0)[lb]{\smash{{\SetFigFont{8}{9.6}{\rmdefault}{\mddefault}{\updefault}{\color[rgb]{0,0,0}$h$}%
}}}}
\put(5551,-3361){\makebox(0,0)[lb]{\smash{{\SetFigFont{8}{9.6}{\rmdefault}{\mddefault}{\updefault}{\color[rgb]{0,0,0}$0$}%
}}}}
\put(6451,-3361){\makebox(0,0)[lb]{\smash{{\SetFigFont{8}{9.6}{\rmdefault}{\mddefault}{\updefault}{\color[rgb]{0,0,0}$1$}%
}}}}
\put(8551,-1786){\makebox(0,0)[lb]{\smash{{\SetFigFont{12}{14.4}{\rmdefault}{\mddefault}{\updefault}{\color[rgb]{0,0,0}$M$}%
}}}}
\put(10801,-1711){\makebox(0,0)[lb]{\smash{{\SetFigFont{12}{14.4}{\rmdefault}{\mddefault}{\updefault}{\color[rgb]{0,0,0}$X_+=h^{-1}([1,+\infty])$}%
}}}}
\put(10801,389){\makebox(0,0)[lb]{\smash{{\SetFigFont{12}{14.4}{\rmdefault}{\mddefault}{\updefault}{\color[rgb]{0,0,0}$h(X\setminus\mathcal{N}):=+\infty$}%
}}}}
\put(17776,389){\makebox(0,0)[lb]{\smash{{\SetFigFont{12}{14.4}{\rmdefault}{\mddefault}{\updefault}{\color[rgb]{0,0,0},}%
}}}}
\put(10801,-586){\makebox(0,0)[lb]{\smash{{\SetFigFont{12}{14.4}{\rmdefault}{\mddefault}{\updefault}{\color[rgb]{0,0,0}$X^o=h^ {-1}(]0,+\infty])$}%
}}}}
\put(10951,-3436){\makebox(0,0)[lb]{\smash{{\SetFigFont{12}{14.4}{\rmdefault}{\mddefault}{\updefault}{\color[rgb]{0,0,0}$M=h^{-1}([0,+\infty])$}%
}}}}
\put(10876,-2536){\makebox(0,0)[lb]{\smash{{\SetFigFont{12}{14.4}{\rmdefault}{\mddefault}{\updefault}{\color[rgb]{0,0,0}$X_-=h^ {-1}]0,1[$}%
}}}}
\end{picture}%

\smallskip A riemannian metric $g$ on $X^{o}$ satisfying:
\begin{equation}\label{product-metric}
 (\varphi_c)_*g(h,y)=dh^2+g_L(y), \ \ (h,y)\in]0,+\infty[\times L
\end{equation} 
on $\cN\setminus\{c\}$ is chosen, where $g_L$ is a riemannian metric on $L$. The corresponding
exponential maps for $X^{o}$ and $L$ are denoted $e$ and $e_L$, and the injectivity radius is assumed to
be greater than $1$ in both cases. The geodesic distances are denoted $\dist,\dist_L$.
The associated riemannian measure on $X^{o}$ and $L$
will be noted $d\mu$ and $d\mu_L$ and the associated
Lebesgue measures on  $T_x X^{o}$ and $T^*_xX^{o}$ for $x\in
X^{o}$ will be noted $dX$ and $d\xi$. We will note $d\mu^\RR$ the Lebesgue measure
on $\RR$.

\smallskip We shall assume that $X^{o}$ is oriented, but all constructions below can be done in
the general case with half densities bundles. 

\smallskip The tangent space of a conical pseudomanifold $X$ was defined in
\cite{DL} by:
 \begin{equation}\label{TcX}
    T^cX=\cC_{X_-}\cup TX_+\rightrightarrows X^{o}
  \end{equation}
The unit space is $X^{o}$. This is a disjoint union where $\cC_{X_-}$ is the pair groupoid of
$X_-$ and $TX_+$ has groupoid structure equal to its vector bundle
structure. We will mainly use in this paper a slightly different but equivalent
(at the level of $K$-theory) definition of the tangent space  which was also
given in \cite{DL}: 
 \begin{equation}\label{TqX}
  T^qX=T]0,1[_h\times \cC_{L}\cup TX_+\rightrightarrows X^{o} 
 \end{equation}  
We will refer to (\ref{TqX}) as the 'q' version and (\ref{TcX}) as the
'c' version of the tangent space of $X$. 

\smallskip The {\sl tangent groupoid} is defined for the 'c' and the 'q'
version by: 
\begin{equation}\label{cGc}
 \cG^c = T^cX\times\{0\}\cup\cC_{X^{o}} \times ]0,1]_t.
\end{equation}
\begin{equation}\label{cGq}
    \cG^q = T^qX\times\{0\} \cup \cC_{X^{o}}\times ]0,1]_t
\end{equation}
In order to write down on $T^qX$ some constructions made in \cite{DL} for
$T^cX$, one needs the following deformation groupoids:
\begin{equation}\label{H} 
   H= T^qX\times\{u=0\} \cup T^c X\times ]0,1]_u,  
\end{equation}
\begin{equation}\label{cH}
 \cH = \cG^q\times\{u=0\} \cup \cG^c\times]0,1]_u
\end{equation}
Let us recall that $\cH$ has three deformation parameters noted
$h,t,u$ and contains all previous groupoids: 
\begin{equation}\label{faces-cH}
  \cH|_{t=0}=H,\ \cH|_{u=0}=\cG^q,\ \cH|_{u=1}=\cG^c,\
  \cH|_{t=0,u=0}=T^qX,\ \cH|_{t=0,u=1}=T^cX.
\end{equation}

\vspace{1cm}

\begin{center}
\includegraphics[width=7cm]{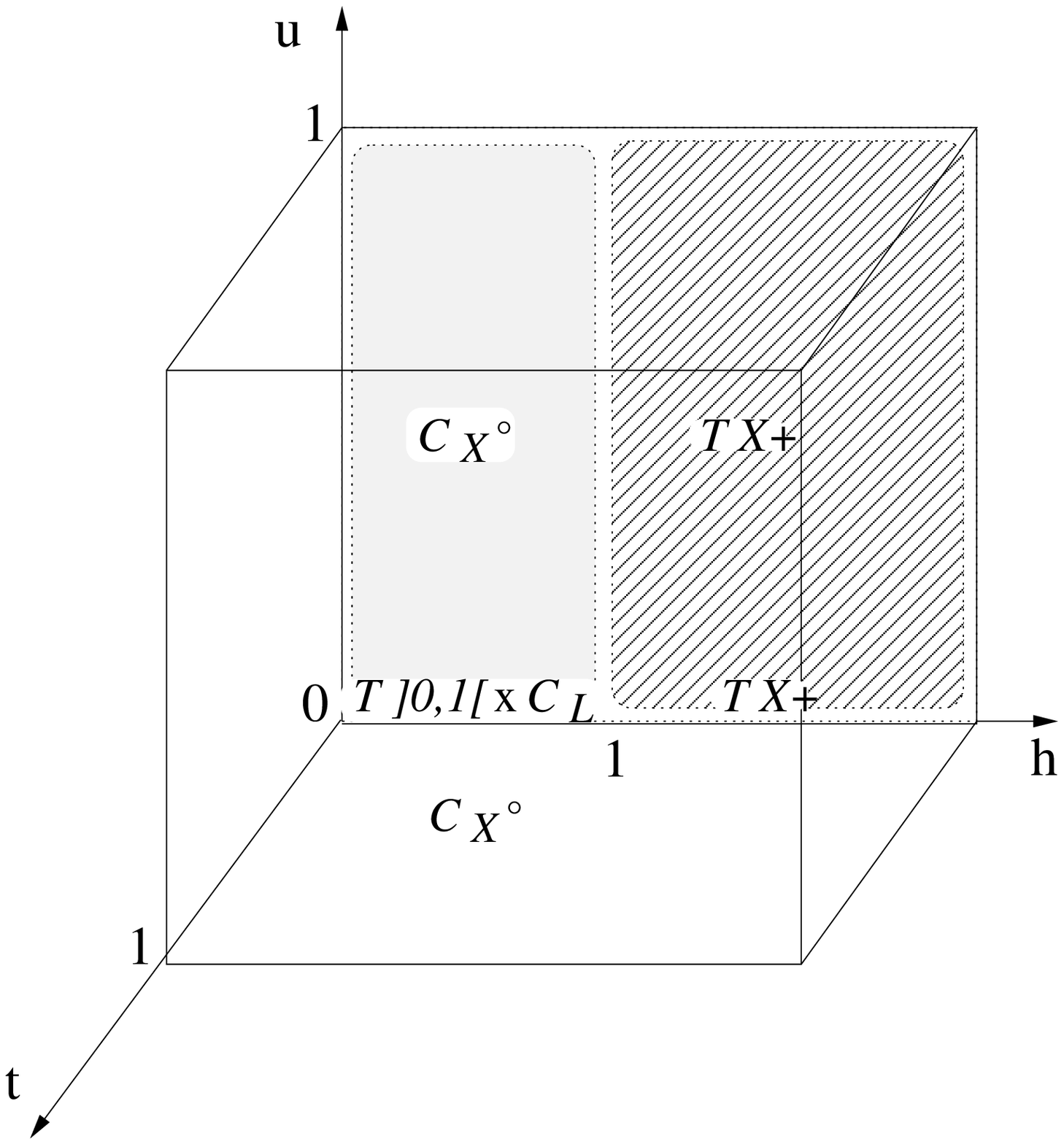}
\end{center}

We describe now a differentiable structure for
$\cH$. This is rather technical, and the unfamiliar reader can skip
this construction up to the remark \ref{exp-cGq} which will be reused later. 

\smallskip The unit space of $\cH$ is $\cH^0=X^{o}\times
[0,1]^2_{t,u}$. We cover $\cH^0$ by four open subsets : 
$\cH^0=A\cup B\cup C\cup D$ with:
\begin{equation}\label{cover-units-of-cH}
 \begin{matrix} 
   A = \interior(X_+)\times[0,1]^2_{t,u} & B =
   X^{o}\times]0,1]_t\times[0,1]_u \\
  C = X_-\times[0,1]^2_{t,u} & D=
  h^{-1}(]1-\varepsilon,1+\varepsilon[) \times[0,\varepsilon[_t\times[0,1]_u
 \end{matrix}
\end{equation}
Here $\varepsilon$ is an arbitrary small number and
$\interior(X_+)=X_+\setminus\partial X_+$. We get the
following cover of $\cH$:
\begin{equation}\label{cover-cH}
 \cH = \cH|_A\cup\cH|_B \cup\cH|_C \cup\cH|_D
\end{equation}
We have $\cH|_A=\cG_{\interior(X_+)}\times[0,1]_u$ where $\cG_{\interior(X_+)}$
is the tangent groupoid \cite{co0} of $\interior(X_+)$. We provide it with
its usual smooth structure and $\cH|_A$ inherits the obvious product
smooth structure. Next, the smooth structure of 
$\cH|_B=\cC_{X^{o}}\times]0,1]_t\times[0,1]_u$ is the product one. For
the two remaining subgroupoids, we need to specify some {\sl gluing}
functions in the deformation parameters. We choose once for all a smooth decreasing function:
\begin{equation}\label{gluingfunction-tau} 
   \tau : [0,+\infty[\to [0,1]
\end{equation}
satisfying $\tau(h)=1$ on $[0,1/2]$
and $\tau^{-1}\{0\}=[1,+\infty[$.
Let:
\begin{equation}\label{gluingfunction-kappa} 
   \kappa : [0,+\infty[\times [0,1]\to [0,1]
\end{equation}
be a smooth function satisfying 
$\min\big(1,\tau(h)+t\big)\le\kappa(h,t)\le 1$ and
 $\kappa(h,t)=\tau(h)+t$ if $\tau(h)+t\le 3/4$.
Let:
\begin{equation}\label{gluingfunction-mu} 
   \mu : [0,+\infty[\times [0,1]\times[0,1]\to [0,1]
\end{equation}
be a smooth function satisfying 
$\min\big(1,u\tau(h)+t\big)\le\mu(h,t,u)\le 1$
and $\mu(h,t,u)=u\tau(h)+t$ if
$u\tau(h)+t\le3/4$.
Let:
\begin{equation}\label{gluingfunction-l} 
   l : \RR_+=]0,+\infty[\to \RR
\end{equation}
be a smooth bijective function satisfying $\frac{d}{d h}l>0$ and $l=\Id$ on a neighborhood of
$[1,+\infty[$ in $]0,+\infty[$.

\smallskip Coming back to the subgroupoid $\cH|_C$, observe that: 
\begin{equation}\label{subgrpd-cHc}
 \cH|_C = \left(T]0,1[_h\times\{(0,0)\}\cup
 \cC_{]0,1[_h}\times[0,1]^2_{t,u}\setminus\{(0,0)\}\right)\times\cC_L
\end{equation}
We note shortly $G_I$ the first factor in $\cH|_C$. Let $\cV$ be the
open subset of $T]0,1[_h\times[0,1]^2_{t,u}$ given by:
 $$
  \cV = \{ (h,\lambda,t,u) \ | \ l^{-1}(l(h)+\mu(h,t,u)\lambda)<1\}
 $$
Then we provide $G_I$ with the smooth structure such that the
bijective map $E_{G_I} : \cV\to G_I$ given by: 
\begin{equation}\label{exp-GI}   
      E_{G_I}(h,\lambda,t,u)=
        \left\lbrace\begin{matrix} \big(h,\lambda,0,0\big) &
         \text{ if } (t,u)=(0,0) \\
         \big(h,l^{-1}(l(h)+\mu(h,t,u)\lambda),t,u\big)&
         \text{ if } (t,u)\not=(0,0) 
       \end{matrix}\right.
\end{equation}
is a diffeomorphism. Thus, $\cH|_C$ inherits the product structure of
$G_I\times\cC_L$. 

\smallskip It is tedious but not difficult to check that the smooth
structures given to $\cH|_A,\cH|_B,\cH|_C$ are compatible on their common
domain and it remains to give $\cH|_D$ with a compatible smooth
structure. 

\smallskip Remember that $X^{o}$ and $L$ are riemmannian with
exponential maps denoted by $e$ and $e^{L}$ (see paragraph \ref{notations}) and
consider now the open subset $\cU$ of $TX^{o}\times[0,1]^2_{t,u}$ 
given by the set of $(x,V,t,u)\in TX^{o}\times[0,1]^2$ satisfying:
\begin{itemize}
  \item if $h(x)\ge1$ then $tV\in\dom(e_{x})$,
  \item if $h(x)<1$ then writing $x=(h,y)\in]0,1[\times L$, 
    $V=(\lambda,W)\in\RR\times T_yL$ under the
    identification $X_-=]0,1[_h\times L$, we 
   have $\kappa(h,t)W\in\dom(e^L_{y})$ and
   $(h,\lambda,t,u)\in\cV$.
\end{itemize} 
We define a  injective map $E_{\cH} : \cU\to \cH$ by setting:
\begin{itemize}
 \item for $(x,V)\in TX^{o}$ with $h(x)\ge 1$:
\begin{equation}\label{exp-cH-h-biggerthan1}
  E_{\cH}(x,V,t,u)=
        \left\lbrace\begin{matrix} \big(x,e_x(tV),t,u\big) &
         \text{ if } t>0 \\
         \big(x,W,0,u\big)&
         \text{ if } t=0 
       \end{matrix}\right.
\end{equation}
\item for $(x,V)\in TX^{o}$ with $h(x)< 1$ and $x=(h,y)$,
$V=(\lambda,W)$ as above:
\begin{equation}\label{exp-cH-h-lessthan1}   
      E_{\cH}(h,y,\lambda,W,t,u)=
        \left\lbrace\begin{matrix} \big(h,y,\lambda,e^L_y(\tau(h)W),0,0\big) &
         \text{ if } (t,u)=(0,0) \\
         \big(h,y,l^{-1}(l(h)+\mu(h,t,u)\lambda),e^L_y(\kappa(h,t)W),t,u\big)&
         \text{ if } (t,u)\not=(0,0) 
       \end{matrix}\right.
\end{equation}
\end{itemize}
One can check that $\cU_D=E_{\cH}^{-1}(\cH|_D)$ is an open subset of
$\cU$. We provide $\cH|_D'=\cH|_D\cap E_{\cH}(\cU)$ with the smooth structure
such that the map : $E_{\cH}:\cU_D\to \cH|_D'$ is a
diffeomorphism. On the other hand,
$\cH|_D''=\cH|_D\cap(\cH|_A\cup\cH|_B\cup\cH|_C)$  is open in 
$\cH|_A\cup\cH|_B\cup\cH|_C$  so it inherits a smooth structure as
a submanifold of $\cH|_A\cup\cH|_B\cup\cH|_C$.  The smooth structure given to 
$\cH|_D'$ and $\cH|_D''$ are compatible and cover $\cH|_D$. The resulting smooth
structure of $\cH|_D$ is compatible with the ones given to the three
other subgroupoids, so we have given to $\cH$ a smooth structure for
which $E_{\cH}$ is an exponential map.  
\begin{rem}\label{exp-cGq} 
\begin{enumerate}
\item The gluing function $l$ is important in the description of the Poincar\'e
  duals of elliptic noncommutative symbols. We will see that choosing $l=\log$ near
  $h=0$ leads to this description with the help of $b$-calculus, while
  choosing $l(h)=-1/h$ near $h=0$ would lead to $c$-calculus.
  Different choices of $l$ produce different but diffeomorphic smooth
  structures on $\cH$. Indeed, if $l,m$ are two such choices, the map $\cH\to \cH$
  equal to identity if $(t,u)\not=(0,0)$ or $h\ge 1$ and sending
  $(h,\lambda,x,y,0,0)$ to $(h,\lambda.\frac{m'(h)}{l'(h)},x,y,0)$ is
  a smooth isomorphism between $\cH$ with the smooth structure given by $l$ and $\cH$
  with the smooth structure given by $m$. This follows from a simple
  but tedious calculation and the fact that for any smooth function
  $f:\RR \to \RR$,  the map: 
  $$\phi : \RR^{3}\to \RR; \ (x, \mu,\lambda)\mapsto  
  \begin{cases}
    \frac{f(x)-f(x-\mu\lambda)}{\mu} & \hbox{ if } \mu\not=0 \\
     \lambda f'(x) & \hbox{ if } \mu=0 
  \end{cases}
$$
is smooth on $\RR^{3}$.  
\item All other gluing functions are technical ingredients and their choice has no incidences on the
  desired description. 
\item All subgroupoids listed in (\ref{faces-cH}) inherit smooth
  structures and exponential maps from those of $\cH$.
\end{enumerate}
\end{rem}
We will often use $\cG^q=\cH|_{u=0}$ in the sequel. Using $E_{\cH}$,
we get an exponential map
for this groupoid:
\begin{itemize}
 \item for $(x,V,t)\in \cU|_{u=0}$ with $h(x)\ge 1$:
\begin{equation}\label{exp-cGq-h-biggerthan1}
  E_{\cG^q}(x,V,t)=
        \left\lbrace\begin{matrix} \big(x,e_x(tV),t,u\big) &
         \text{ if } t>0 \\
         \big(x,W,0,u\big)&
         \text{ if } t=0 
       \end{matrix}\right.
\end{equation}
\item for $(x,V,t)\in \cU|_{u=0}$ with $h(x)< 1$ and $x=(h,y)$,
$V=(\lambda,W)$:
\begin{equation}\label{exp-cGq-h-lessthan1}   
      E_{\cG^q}(h,y,\lambda,W,t)=
        \left\lbrace\begin{matrix} \big(h,y,\lambda,e^L_y(\tau(h)W),0,0\big) &
         \text{ if } t=0 \\
         \big(h,y,l^{-1}(l(h)+t\lambda),e^L_y(\kappa(h,t)W),t,u\big)&
         \text{ if } t>0
       \end{matrix}\right.
\end{equation}
where we have replaced $\mu(h,t,0)$ by $t$ to simplify. 
\end{itemize}
The inverse of the exponential map
$E_{\cG^q}$ will be noted shortly $\Theta$.  

\smallskip The following define a Haar
system for $\cG^q$ which is necessary to define in a
convenient way its $C^*$-algebra:
\begin{equation}\label{Haar-system-cGq} 
   \begin{array}{ccc}  
    t>0, & \cG^q_{(x,t)}=X^{o}\times\{t\}, &  
       d\lambda^{(x,t)}(x')=\frac{l'(h')}{t\tau(t,h')^n}d\mu_{x'}
       =d\lambda^t  \ (h'=h(x'))\\
       t=0 \text{ and } h<1, & \cG^q_{(h,y,0)}=\RR\times L\times\{0\},  & 
      d\lambda^{(h,y,0)}=\frac{1}{\tau(h)^n}d\mu^\RR d\mu^L =d\lambda^{h,0}\\  
    t=0 \text{ and } h\ge1, & \cG^q_{(x,0)}=T_xX^{o}, & d\lambda^{(x,0)}=d\nu_{x}=d\lambda^x 
   \end{array}   
\end{equation}
Remark that $d\lambda^1$ is equal to $\frac{1}{h}d\mu$ near $h=0$, in
other words it coincides with the density coming from a
$b$-metric like: 
 \begin{equation}\label{bmetric}
   g_b(h,y) = \frac{dh^2}{h^2} + g_L(y)
 \end{equation}

\section{Equivalences of tangent spaces and Dirac elements}\label{section2}
\subsection{Two equivalent tangent spaces} 
The main results of \cite{DL} are the construction 
of a Dirac element $D^c\in KK(T^cX\times X,\cdot)$, a dual
Dirac element $\lambda^c\in KK(\cdot,T^cX\times X)$, where $\cdot$
stands for a point space, and the 
computations in bivariant $K$-theory:
\begin{equation} 
    \lambda^c\underset{T^cX}{\otimes} D^c = 1\in KK(X,X)
 \text{ and }
 \lambda^c\underset{X}{\otimes} D^c = 1\in KK(T^cX,T^cX)
\end{equation}
which give in particular the isomorphism
 $$
  \Sigma^c=(\lambda^c\underset{X}{\otimes}\cdot)=(.\underset{T^cX}{\otimes} D^c)^{-1}
 $$
in (\ref{PD-map-conic-case}). In this paper, we will prefer to work
with $T^qX$ rather than $T^cX$, because the analog for $T^qX$ of $\Sigma^c$:
\begin{equation}\label{PD-map-conic-q-case} 
  \Sigma^q=(\lambda^q\underset{X}{\otimes}\cdot)=(.\underset{T^qX}{\otimes} D^q)^{-1}: \ 
  K_0(X)\overset{\simeq}{\longrightarrow} K^0(T^qX)
 \end{equation} 
has a nice description. We are going to describe the $KK$-equivalence
between $T^qX$ and $T^cX$ in order to have a correct representant of
$D^q$. 
\begin{prop}\label{KK-equiv-qc} The deformation element $\partial_H\in
  KK(T^qX,T^cX)$ associated with  $H$ in
  (\ref{H}) is a $KK$-equivalence. 
\end{prop}
\begin{pf}
Let $\ev[+] : H\to TX_+\times[0,1]$ be the
restriction morphism and consider the commutative diagram:
\begin{equation}\label{CD-H}
 \begin{CD}  
    0 & & 0 & & 0 \\
    @AAA  @AAA  @AAA \\
  C^*(TX_+) @<\ev[u=1]<< C^*([0,1]_u\times TX_+) 
  @>\ev[u=0]>> C^*(TX_+) \\
  @A\ev[+]AA  @A\ev[+]AA @A\ev[+]AA \\
   C^*(T^cX)  @<\ev[u=1]<< C^*(H)  @>\ev[u=0]>> C^*(T^qX) \\
     @AAA  @AAA  @AAA \\
  C^*(\cC_{]0,1[\times L}) @<\ev[u=1]<< C^*(H|_{h<1}) @>\ev[u=0]>> C^*(T]0,1[\times \cC_L) \\
  @AAA  @AAA  @AAA \\
 0 & & 0 & & 0
 \end{CD}
\end{equation}
The columns are exact. On the top line the induced maps in $K$-theory provide 
\begin{equation}
  \label{topinvertibility}
[\ev[u=0]]^{-1}\ot[\ev[u=1]]=1\in KK(TX_+,TX_+)  
\end{equation}
Observe that
\begin{equation}
  \label{Hmoins}
   H|_{h<1}= T]0,1[\times\cC_L \times\{0\} \cup \cC_{]0,1[\times L}\times
   ]0,1] \simeq \left(T]0,1[\times\{0\} \cup \cC_{]0,1[}\times]0,1]_u\right)\times\cC_L,
\end{equation}
and that $\cG_{]0,1[}=T]0,1[\times\{0\} \cup \cC_{]0,1[}\times]0,1]$
is the usual tangent groupoid of the manifold $]0,1[$. The associated
Kasparov element
 $$
   \partial_{\cG_{]0,1[}}\in KK(T]0,1[,\cC_{]0,1[})\simeq
   KK(C_0(\RR^2),\CC)\simeq \ZZ
 $$
is invertible with inverse given by the Bott generator of
$KK(\CC,C_0(\RR^2))$. It follows that in the bottom line 
\begin{equation}
  \label{bottominvertibilty}
   [\ev[u=0]]^{-1}\ot[\ev[u=1]]=\mathbf{s}_{C^*(\cC_{L})}\left(\partial_{\cG_{]0,1[}}\right)\in 
   KK(T]0,1[\times\cC_L,\cC_{]0,1[}\times\cC_L)
\end{equation}
is invertible. Here $\mathbf{s}_A:KK(B,C)\to KK(B\ot A,C\ot A)$ is the usual
tensorisation operation in Kasparov theory. 

\smallskip In particular the Kasparov elements $[\ev[u=1]]$ in
(\ref{topinvertibility}) and (\ref{bottominvertibilty}) are
invertible. Hence, applying the five lemma to the long exact sequences
in $KK$-theory associated with the
first two columns of (\ref{CD-H}) give the invertibility of the
element $[\ev[u=1]]$ in the middle line of (\ref{CD-H}). This yields
that $\partial_H$, equal to $[\ev[u=0]]^{-1}\ot[\ev[u=1]]$ in the
middle line, is invertible. 
\end{pf}

\subsection{The Dirac element for the 'q' version}
Let us turn to the description of $D^q$. We define:  
 \begin{equation}\label{df-Dq}
   D^q:=\partial_H\underset{T^cX}{\ot} D^c \in KK(T^qX\times X, \cdot)
\end{equation}
We recall the construction of $D^c$. We set: 
 $$ 
  \partial^c:= \partial_{\cG^c}\ot \nu \in KK(T^cX,\cdot) 
 $$
where $\nu$ is the Morita equivalence $\nu\in KK(\cC_{X^{o}},\cdot)$ given by:
 $$
   \nu = \left(L^2(X^{o}),m,0\right)\in KK(\cC_{X^{o}},\cdot)
 $$
$\partial^c$ is called the pre-Dirac element and the Dirac element is
 $$ 
   D^c:= \Phi^c\otimes \partial^c \in KK(T^cX\times X,\cdot)
 $$
where $\Phi^c$ is the $KK$-element associated to the map defined by
\begin{equation}\label{Psic}
 \begin{matrix} 
  \Phi^c\ : &  C^*(T^cX)\ot C(X) & \longrightarrow & C^*(T^cX) \\
     &  a\ot f & \longmapsto & a.f\circ\pi^c
 \end{matrix}
\end{equation}
and $\pi^c$ is the composition of the range map
of $T^cX$ with the projection $X^{o}\to X^{o}/\overline{X_-}\simeq X$. In the 'q' version we set: 
 $$
   \partial^q=\partial_{\cG^q}\otimes\nu
 $$
and
\begin{equation}\label{Psiq}
 \begin{matrix} 
  \Phi^q\ : &  C^*(T^qX)\ot C(X) & \longrightarrow & C^*(T^qX) \\
     &  a\ot f & \longmapsto & a.f\circ\pi^q
 \end{matrix}
\end{equation}
where $\pi^q$  projects $T^qX$ onto $X$ like $\pi^c$ does $T^cX$ on $X$. 
We check that : 
\begin{prop} The following equality holds : 
 $$ 
  D^q= \Phi^q\otimes \partial^q \in KK(T^qX\times X,\cdot) 
 $$
where $D^q$ is defined by (\ref{df-Dq}). 
 \end{prop}
\begin{pf}
Let us consider the commutative diagram:
\begin{equation}\label{CD-cqH} 
\begin{CD}
  C^*(\cC_{X^{o}}) @<\ev[t=1]<<  C^*(\cG^q) @>\ev[t=0]>> C^*(T^qX) \\
  @A\ev[u=0]AA @A\ev[u=0]AA @A\ev[u=0]AA \\
   C^*(\cC_{X^{o}}\times\lbrack 0,1\rbrack) @<\ev[t=1]<< C^*(\cH)
   @>\ev[t=0]>> C^*(H) \\ 
    @V\ev[u=1]VV @V\ev[u=1]VV @V\ev[u=1]VV \\
  C^*(\cC_{X^{o}}) @<\ev[t=1]<< C^*(\cG^c) @>\ev[t=0]>>   C^*(T^cX)
\end{CD}
\end{equation} 
At the level of $KK$-theory, the bottom line of the diagram gives
$\partial^q$, up to the Morita equivalence $\nu$, while the top
line gives $\partial^c$ (up to $\nu$). The right column gives the
$KK$-equivalence $\partial_H$ while the product
$\ev[u=0]^{-1}\ot\ev[u=1]$ in the left column is obviously $1$ in
$KK(\cC_{X^{o}},\cC_{X^{o}})$.  
This gives:
 $$ 
    \partial_H\ot\partial^c=\partial^q.
 $$
To finish, let us introduce the multiplication morphism:
 $$
     \Phi^H : C^*(H)\ot C(X)\To C^*(H) 
  $$ 
given by $\Phi^H(a,f)(\gamma)=a(\gamma)f(\pi^H(\gamma))$, where 
the projection map $\pi^H : H\To X$ extends $\pi^q$ and $\pi^c$ in the
obvious way. Denoting in the same way the restriction morphisms for
the product groupoid $H\times X$, we
get the formulas: 
 $$ 
    \mathbf{s}_{C(X)}\partial_H= \mathbf{s}_{C(X)} ([\ev[u=0]]^{-1}\ot
    [\ev[u=1]])=[\ev[u=0]]^{-1}\ot [\ev[u=1]],
  $$
  $$
    [\ev[u=1]]\ot \Phi^c = \Phi^H \ot [\ev[u=1]]\quad ; \quad
    [\ev[u=0]]^{-1} \ot \Phi^H = \Phi^q \ot [\ev[u=0]]^{-1}.
  $$ 
Hence:
\begin{eqnarray}
  D^q=\partial_H\underset{T^cX}{\ot} D^c &=& \mathbf{s}_{C(X)}(\partial_H) \ot D^c 
    = [\ev[u=0]]^{-1}\ot [\ev[u=1]]\ot \Phi^c\ot
   \partial^c_X\\
   \nonumber 
   &=& \Phi^q \ot[\ev[u=0]]^{-1}\ot[\ev[u=1]] \ot\partial^c
   = \Phi^q\ot\partial^q
\end{eqnarray}
\end{pf}

\section{Cycles of the $K$-theory of the tangent space} 
\subsection{Symbols on a manifold as operators on the tangent space}\label{symbols-manifolds}

To proceed, we need some definitions about pseudodifferential calculus on
groupoids. The notions summed up below can be found with full details in the litterature: \cite{Vas,MP,NWX,CoS}. 

\smallskip Let $G$ be a smooth groupoid (the space of units is
allowed to be a manifold with boundary,  but the fibers are manifolds
without boundary).
 Let $U_{\gamma} : C^{\infty}(G_{s(\gamma)})\To
C^{\infty}(G_{r(\gamma)})$ be the isomorphism induced by right
multiplication:  $U_{\gamma}f(\gamma')=f(\gamma'\gamma)$. A linear operator
$P : C^{\infty}_{c}(G) \to C^{\infty}(G)$ is a {\it $G$-operator} if there exists a family  
$P_{x} : C^{\infty}_{c}(G_{x}) \to C^{\infty}(G_{x}) $ such that
$P(f)(\gamma)=P_{s(\gamma)}(f\vert_{G_{s(\gamma)}})(\gamma)$ and
$U_{\gamma}P_{s(\gamma)}=P_{r(\gamma)}U_{\gamma}$.  

\smallskip A $G$-operator $P$ is a {\it pseudodifferential operator on $G$} (resp. of
order $m$) if for any open local chart  
$\Phi :\Omega \To s(\Omega)\times W$  of $G$ such that 
$s=pr_1\circ \Phi$ (that is, for any
distinguished chart) and any cut-off
function $\chi \in C^\infty_c(\Omega)$,  we have
$(\Phi^{*})^{-1}(\chi P \chi)_{x}\Phi^{*}= a(x,w,D_w)$  where 
$a(x,w,\xi)\in S^*(s(\Omega)\times T^*W)$ is a classical symbol 
(resp. of order $m$).

\smallskip One says that $P$ has {\it support} in $K\subset G$  if 
$\mbox{supp}(Pf)\subset K.\mbox{supp}(f)$ for all $f\in
C^{\infty}_{c}(G)$.   

\smallskip These definitions extend immediately to the case of operators acting
between sections of bundles on $G^{(0)}$ pulled back to $G$ with the
range map $r$. The space of compactly supported pseudodifferential
operators on 
$G$ acting on sections of $r^*E$ and taking values in sections of $r^*F$
will be noted $\Psi_c^*(G,E, F)$. If $F=E$ we get an algebra
denoted by $\Psi_c^*(G,E)$. 

\smallskip Basic examples of the usefulness of these operators are the case of
foliations \cite{CoS,Vas} and manifolds with corners \cite{Mo2}.
This calculus is also used in \cite{DL} to define $KK$-theory classes
and to compute some Kasparov products. Here, 
to motivate our definition of {\it noncommutative symbols} on a singular 
manifold, we explain in more details what is suggested in the
introduction.

\smallskip Let $V$ be a smooth compact riemannian manifold, $E$ a smooth vector
bundle over $V$ and consider the tangent space $TV$ as a smooth
groupoid (thus $r$ and $s$ are equal to the canonical projection map 
$TV\to V$). 

\smallskip  Let $a\in\Psi_c^*(TV,E)$.  By definition,
$a$ is a smooth
family $(a_{x})_{x\in V}$ where $a_{x}$ is a
translation invariant pseudodifferential operators on $T_xV$
(with coefficients in $\End E_x$) and thus can be regarded as a
distribution $a_{x}(X)$ on $T_{x}V$ acting by convolution on
$C^{\infty}_{c}(T_{x}V, E_{x})$, so:
 $$ u\in C^{\infty}_{c}(TV, E), \quad a(u)(x, X)=a*u(x, X)=\int_{Y\in
 T_{x}V} a(x, X-Y)u(Y)dY$$
where the last integral is understood in the distributional
sense. The distribution $a_{x}(X)$ being compactly supported, it has
a Fourier transform $\widehat{a_{x}}(\xi)$ which is just its symbol. The
whole family $(\widehat{a_{x}})_{x}$ identifies with a classical
symbol on $V$ taking values in $\End E)$,  that is,  $\widehat{a}\in
S^*(V,\End E)$ and since the Fourier transform exchanges convolution with
pointwise multiplication,  we get an algebra homomorphism:  
 \begin{equation}\label{symbolsV}
   \begin{array}{cccc}
   \cF : & \Psi_c^*(TV,E) &\To & S^*(V,\End E) \\
         &  a  &\longmapsto & \widehat{a}(x,\xi)
   \end{array}
 \end{equation}
which is obviously injective. Conversely, the inverse Fourier
transform associates to any symbol $b(x, \xi)\in S^*(V,\End E) $ a
distribution $\overset{\vee}{b}(x, X)$ which, as a convolution
operator, is a $TV$-pseudodifferential operator by the formula:
 $$u\in C^{\infty}_{c}(TV, E), \quad \overset{\vee}{b}*u(x, X)=
 \int_{T_{x}V\times T^{*}_{x}V}e^{i(X-Y). \xi} b(x,
 \xi)u(x,Y)dYd\xi $$
Moreover,  introducing a smooth function $\phi(x, X)$ on $TV$ equal to
$1$ if $X=0$ and equal to $0$ if $|X|>1$,  we get:
 $$ 
  \overset{\vee}{a} = \phi. \overset{\vee}{a}+(1-\phi)\overset{\vee}{a}\in   \Psi_c^*(TV,E) + \cS(TV,\End E)
 $$
where $\cS(TV, \End E)$ stands for the (Schwartz) space of smooth
sections whose partial derivatives are rapidly decaying in the
fibers of $TV$.  



\smallskip Thus, enlarging $\Psi_c^*(TV,E)$ as follows:
    $$ \Psi^*(TV,E) := \Psi_c^*(TV,E) + \cS(TV,\End E),  $$
the algebra monomorphism (\ref{symbolsV}) extends to an
algebra isomorphism
\begin{equation}\label{symbolsV+}
     \cF \ :\ \Psi^*(TV,E)\longrightarrow S^*(V,\End E)
 \end{equation}
which preserves the filtrations. For a general discussion about the
enlargement of spaces of compactly supported pseudodifferential
operators by adding regularizing ones,  see \cite{Vas,LMN2005}. 
             

 One can 
then reformulate the classical description of the $K$-theory
with compact supports of $T^*V$:
\begin{prop}\label{KTVbySymbols} 
  Every element in $K^0(T^*V)\simeq K^0(TV)$ has a representant of
  the following form:
   $$
     [a] =  \left( C^*(TV,E\oplus F), \begin{pmatrix}0 & a \\ b & 0
     \end{pmatrix}\right)\in KK(\cdot,TV) 
   $$
   where $E,F$ are smooth vector bundles over $V$ and
   $a\in\Psi^0(TV,E,F)$, $b\in\Psi^0(TV,F,E)$ satisfy
   $ba-1\in\Psi^{-1}(TV,E)$, $ab-1\in\Psi^{-1}(TV,F)$. 
\end{prop}

\subsection{Noncommutative symbols and their ellipticity on a conical pseudomanifold}\label{symbols-pseudomanifolds}
Motivated by the previous approach,  we enlarge the
space of compactly supported pseudodifferential operators on $T^qX$
and define them as noncommutative symbols on the pseudomanifold $X$. 
Definitions are given in the scalar 
case since the presence of vector bundles bring no issues.
We introduce:
\begin{equation}\label{TqXbar}
    \Tb = \{0\}\times\RR\times\cC_L \cup T^qX =T[0,1[_h\times\cC_L\cup
    TX_+ \rightrightarrows M=\overline{X^o}
\end{equation}
\begin{defn}\label{def.schwartz.TqXbar}
 Let $\tau$ be the function choosed in (\ref{gluingfunction-tau})
 and define the function $|.| : T^qX\to \RR_+$ by: 
 $$|\gamma| =\begin{cases}
   \sqrt{\left(\fraq{\mbox{dist}_L(x,y)}{\tau(h)}\right)^2+\lambda^2}& \hbox{ if }
   \gamma=(h,\lambda,x,y)\in T]0,1[\times\cC_L \ (ie, \ h<1) \\
     \sqrt{g_x(X,X)} & \hbox{ if } \gamma=(x,X)\in TX_+ \ (ie, \ h\ge 1)
       \end{cases}
 $$
\end{defn}
The restriction at $t=0$ of the Haar
system $\cG^{q}$ defined in (\ref{Haar-system-cGq})  provides a Haar
system for $T^{q}X$, and extended at $h=0$ in the
obvious way, we get a Haar system for $\Tb$. It is then easy to check
that $|. |$ is a length function with polynomial growth on $\Tb$
 and the corresponding Shwartz algebra is denoted by
$\cS(\Tb)$ (\cite{LMN2005}).  Using the seminorms:
 $$ p_{D,N}(f) = \sup_{\gamma\in T^qX}(1+|\gamma|)^{N}|Df(\gamma)|, 
 $$
where $N$ is a positive integer and $D\in\hbox{Diff}(\Tb)$
is a differential operator on $\Tb$, this Schwartz algebra can be
presented as follows: 
$$ \cS(\Tb) = \{ f\in C^\infty(\Tb) \ | \ p_{D,N}(f)<+\infty\ \forall
N\in \NN, \ \forall D\in \hbox{Diff}(\Tb)\}
$$ 

By restriction at $h=0$ of these functions, we get the Schwartz
algebra $\cS(\RR\times\cC_L)$ of the groupoid $\RR\times\cC_{L}$
endowed with the restricted Haar system. 


\begin{defn}\label{symbol} The algebra of noncommutative symbols on $X$ is defined by: 
\begin{equation}\label{def-symb-alg}
 S^*(X) = \Psi_c^*(\Tb) + \cS(\Tb) \subset \Psi^*(\Tb)
\end{equation}
We also define $S^*_0(X)$ as the kernel of the restriction
homomorphism at $h=0$:
\begin{equation}\label{restrictions-symb} 
    \begin{matrix}
       \rho \ : &  S^*(X) & \To & \Psi^*_c(\RR\times\cC_L)+\cS(\RR\times\cC_L)\\
                &   a     & \longmapsto \ a|_{h=0} 
                    \end{matrix}  
 \end{equation}
\end{defn}
\begin{rem}
\item[{\bf (a)}] The image of $\rho$ is exactly the algebra
  $\cP^*_{inv}(\RR\times L)$
of translation invariant pseudodifferential operators on $\RR\times L$
defined by R. Melrose in \cite{Mel2}.

 \item[{\bf (b)}] The smoothness of noncommutative symbols up to $h=0$ can be relaxed and
    singular behaviors can be of interest, see paragraph
   \ref{FuchsSymbols}.
\end{rem}
The following observations will lead to the notion of ellipticity for
these noncommutative symbols.
\begin{prop}\label{symbols-as-multipliers}
The following inclusion holds
\begin{equation}\label{symbols-are-multipliers}
     S^0(X)\subset \cM(C^*(T^qX)).
\end{equation}
Moreover, we have $S^{-1}(X)\subset C^*(\Tb)$ and
\begin{equation}\label{symb-1-0-are-compact}
   S^{-1}_0(X) = S^{-1}(X)\cap C^*(T^qX).
\end{equation}
\end{prop}
\begin{pf} It is known from \cite{MP,Vas} that
  $\Psi^0(\Tb)\subset\cM(C^*(\Tb))$ and 
   $\Psi^{-1}(\Tb)\subset C^*(\Tb)$. Since
   $C^*(T^qX)$ is an ideal of $C^*(\Tb)$ and  
  since any $a\in \Psi^0(\Tb)$ maps  $C^*(T^qX)$ to itself, (\ref{symbols-are-multipliers}) is true.
Since $C^*(T^qX)$ is the kernel of the restriction morphism
$C^*(\Tb)\to C^*(\RR\times\cC_L)$ at $h=0$, (\ref{symb-1-0-are-compact}) is obvious. 
\end{pf} 

In the sequel, the algebra of small $b$-calculus 
\cite{Mel} will be denoted by $\cP^*_b(M)$ and its ideal consisting of operators with vanishing
indicial families will be denoted by $\cP^*_{b,0}(M)$. Given a
$\cG^{q}$-pseudodifferential operator $P$, its restriction $P|_{t}$
at any $t>0$ is a $\cC_{X^{o}}$-pseudodifferential operator,  that
is, an ordinary pseudodifferential operator on the (open) manifold $X^{o}$. 
In fact,  we will denote by $\Psi^*_b(\cG^q)$  the algebra of
$\cG^{q}$-pseudodifferential operators whose restrictions $P|_{t}$ at any
$t>0$ are in the $b$-calculus of $M=\overline{X^{o}}$,  that is,
such that for all $t>0$,  $P|_{t}\in\cP^*_b(M)$. 
The ideal of operators $P\in\Psi^*_{b}(\cG^q)$ such that $P|_{t}\in\cP^*_{b,0}(M)$ for
all $t>0$ will be denoted by $\Psi^*_{b,0}(\cG^q)$. The previous proposition
extends to these spaces: 
\begin{prop}\label{cGq-operators-as-multipliers}
The following inclusions hold: 
\begin{equation}\label{cGqpsib-are-multipliers} 
    \Psi^0_{b}(\cG^q)\subset\cM(C^*(\cG^q))
\end{equation}
\begin{equation}\label{cGqpsib0-are-compact} 
    \Psi^{-1}_{b,0}(\cG^q)\subset C^*(\cG^q)
\end{equation}
\end{prop}
This follows from properties of $b$-calculus and the proposition \ref{symbols-as-multipliers}.
\begin{defn}
 A noncommutative symbol $a\in S^*(X)$ is elliptic if it is invertible in  $S^*(X)$
 modulo $S^{-1}_0(X)$. \\
 A noncommutative symbol $a\in S^*(X)$ is relatively elliptic if it is invertible in
 $S^*(X)$  modulo $S^{-1}(X)$.
\end{defn}
The relative ellipticity of $a\in S^*(X)$ is exactly its ellipticity
as a pseudodifferential operator on $\Tb$. The notion of ellipticity
for our noncommutative symbols is stronger and is similar to the notion
of full ellipticity \cite{Mel,Mo2003}.

\smallskip Indeed, let $\sigma(a)\in C^\infty(S^*M)$ be the principal symbol of
$a\in S^*(X)$ viewed as a pseudodifferential operator on $\Tb$. We call 
$(\sigma(a),\rho(a))\in C^\infty(S^*M)\times\cP^*_{inv}(\RR\times L)$ 
 the {\sl leading part of} $a$.  
\begin{prop} The following assertions are equivalent:
\begin{enumerate}
\item The noncommutative symbol $a$ is elliptic on $X$.
\item  The leading part of $a$ is invertible.
\end{enumerate}
\end{prop}
 \begin{pf} (i)$\Rightarrow$ (ii) is obvious. Conversely,
let $a\in \Psi^d(\Tb)$ be a noncommutative symbol whose principal part
 is invertible. Since $a$ is an elliptic
$\Tb$-pseudodifferential operator, we can choose $\widetilde{b}\in
\Psi^{-d}(\Tb)$ inverting $a$ modulo $\Psi^{-1}(\Tb)$.  From the
smoothness of the family $a$ 
we get a continuous map 
$h\in[0,1]\mapsto a|_h\in \cP^*_{inv}(\RR\times L)$ and the invertibility of $\rho(a)$ implies the invertibility of
$a|_h$ if $h<\alpha$ for some $\alpha>0$. We pick a cut-off function
$\omega \in C^\infty_c[0,\alpha[$ such that $\omega(0)=1$  and we set :
 $$ 
   b = \omega (a|_{h})^{-1}+(1-\omega)\widetilde{b}.
 $$
Then
  $$
     ab = \omega + (1-\omega)a\widetilde{b} =
       \omega +(1-\omega)(1+q) = 1 +(1-\omega)q,
  $$ 
where $q\in\Psi^{-1}(\Tb)$, and  $ab-1\in S^{-1}_0(X)$ is proved. Things
are similar for $ba-1$. 
\end{pf} 

\subsection{$K$-theory of the tangent space}

We prove in this paragraph that elliptic noncommutative symbols, when vector bundles
are allowed, are the cycles of  $K^0(T^qX)$. As already quoted,
definitions \ref{def.schwartz.TqXbar} and \ref{symbol} extend
immediately to the case of vector bundles and we note $S^*(X,E,F)$ the
space of noncommutative symbols on $X$ acting between sections of bundles $E,F$ over
$M$. With the convention $S^*(X,E):=S^*(X,E,E)$, the proposition
\ref{symbols-as-multipliers} becomes: 
 $$ 
   S^0(X,E)\subset\cL(C^*(T^qX,E)),\qquad
   S^{-1}_0(X,E)\subset\cK(C^*(T^qX,E)) 
 $$
Therefore, we can associate to each elliptic noncommutative symbol 
$a\in S^0(X;E,F)$ on $X$ of order $0$, an element in the
$K$-theory of $T^qX$: 
\begin{defn}\label{bounded-symbol} 
Let $a\in S^0(X;E,F)$ be an elliptic noncommutative symbol on $X$. We set:
$$[a]:=[C^*(T^qX,E\oplus F),\mathbf{a}]\in KK(\cdot,T^qX)\simeq K^0(T^qX)$$
where:  
$$\mathbf{a}=\begin{pmatrix} 0 & b\\ a & 
    0\end{pmatrix}$$ 
and $b$ is any noncommutative symbol inverting  $a$ modulo $S^{-1}_0$.
\end{defn}
It is straightforward that $[a]$ does not 
depend on the choice of the quasi-inverse $b$.  
The main result of this section is that proposition \ref{KTVbySymbols}
holds in this new framework:   
\begin{thm}\label{symbgenerate} Every element of $K^0(T^qX)$ has a
  representant among elliptic noncommutative symbols. More precisely:
 $$ K^0(T^qX) = \{ [a] \ | \ a \text{ is an elliptic noncommutative symbol on } X
\text{of order } 0\},
 $$
\end{thm}
\begin{rems} 
\begin{enumerate}
  \item Considering the Kasparov ungraded modules given by 
 $ (C^*(T^qX,E),a)$ where $a\in S^0(X,E)$ and
   $a^2-1\in S^{-1}_0(X,E)$, the conclusion is the same for $K^1(T^qX)$.
  \item In the same way, relative elliptic noncommutative symbols span the $K$-theory
    of $\Tb$. Observe that $\Tb$ is $KK$-equivalent to $TX^{o}$ which is
    $KK$-dual to $M=\overline{X^{o}}$. 
\end{enumerate}
\end{rems} 
{\bf Proof of the theorem: }
Let us denote by $\Delta$ the following subset of $K^0(T^qX)$:
$$\{ [a] \ | \ a \text{ is an elliptic noncommutative symbol on } X
\text{of order } 0\}. $$
From the exact sequence of $C^*$-algebras:
\begin{equation}\label{ES-TqX-TX_+} 
   0\to C^*(\cC_L\times TI)\overset{i}{\to}
   C^*(T^qX)\overset{\ev[+]}{\to}  
   C^*(TX_+)\to  0,
\end{equation}
we get the exactness of 
 $$ 
  K^0(\cC_L\times TI)\overset{i}{\To}
  K^0(T^qX)\overset{\ev[+]}{\to} K^0(TX_+). 
 $$
If $i(K^0(\cC_L\times TI))\subset\Delta$ and
$\ev[+](K^0(T^qX))\subset\ev[+](\Delta)$ then the theorem is true. 
These inclusions are checked in the following lemmas.  

\begin{lem}\label{im-i-ev+}
The inclusion $i(K^0(\cC_L\times TI))\subset\Delta$ holds.
\end{lem}
{\bf Proof of the lemma: }
It is sufficient to find a generator $e$ of $K^0(C^*(\cC_L\times
TI))\simeq\ZZ$ such that $i_*(e)\in\Delta$. We will define first an
appropriate generator of $K_0(C_0(\RR^2)\ot \cK(L^2(L)))$ and then we
will use  an isomorphism $C^*(\cC_L\times TI)\simeq C_0(\RR^2)\ot \cK(L^2(L))$,
 
\smallskip Let us choose the following generator of $K_0(C_0(\RR^2))$: 
 $$
   x = [\cE, 1, F] \hbox{ where: } \cE = C_0(\RR^2)\oplus C_0(\RR^2),
   \ F := \frac{d}{\sqrt{1+d^2}} \hbox{ and } d := 
  \begin{pmatrix} 0 & h-i\lambda \\ h+i\lambda & 0 \end{pmatrix}.
 $$  
On the other hand, let $B_+$ be an elliptic pseudodifferential
operator on $L$ with index $1$. Without loss of generality, we may assume that $B_+$ is
of order $1$, is almost unitary (ie, unitary modulo $0$ order
operators) and acts between sections of a trivial bundle 
$L\times \CC^k$. Let $b_+\in C^\infty(S^*L, U_k(\CC))$ be its principal
symbol. Then the following represents $1\in K^0(\cC_L)$: 
 $$
   x' = [\cE', 1, F'],\hbox{ where: } \cE' =
   \cK((L^2(L;\CC^k))^2) , \ 
   F' := \frac{B}{\sqrt{1+B^2}} \hbox{ and } B := 
   \begin{pmatrix} 0 & B_+^* \\B_+ & 0 \end{pmatrix}; 
 $$ 
Now the Kasparov product $x''=x\underset{\CC}\ot x'$ is a generator 
of $K_0(C_0(\RR^2)\ot\cK(L^2(L)))$  and is represented by 
$(\cE",F")$ where: 
 $$
  \cE''= \cE\underset{C_0(\RR^2)}{\hot}(C_0(\RR^2)\ot\cE') \simeq 
  \cE\underset{\CC}\ot \cE', \qquad F" = \frac{D}{\sqrt{1+D^2}} 
  \quad \text{ and } D = d\hot I_{2k}+I_2\hot B.
 $$ 
Here $I_n$ denotes the identity matrix of rank $n$ and $\hot$ is the graded
tensor product. Recall the matricial expression of
$D$: 
 $$
    D= d\hot I_{2k}+I_2\hot B= 
   \begin{pmatrix} 0 & 0 & 1\ot B_+^{*} & d_-\ot 1 \\
                             0 & 0 & d_+\ot 1 & -1\ot B_+ \\
                             1\ot B_+ & d_-\ot 1 & 0 &0  \\
                             d_+\ot 1 & -1\ot B_+^{*} & 0 & 0 
   \end{pmatrix} 
 $$
It is clear that $D$ is a pseudodifferential operator on $L$ with
parameters $(h,\lambda)\in \RR^2$ (of order $1$) in the sense of \cite{Shub}, acting on the sections of the
product bundle $L\times\CC^{4k}$. Following the construction of complex powers given in
\cite{Shub}, we see that $F"$ remains in the same space of operators
with parameters (but of course, it is of order $0$). 
Let us find a better representant of $x''$  by trivializing $F''$ at $+\infty$. Let us introduce the matrix:
 $$
    J=\begin{pmatrix} 0 & K \\ 
   K & 0 \end{pmatrix}\in M_{4k}(\CC) \text{ where } 
   K = \begin{pmatrix} 0 &  I_k \\ 
                     I_k & 0 
       \end{pmatrix}\in M_{2k}(\CC).
 $$ 
We choose a smooth decreasing function $M$ equal to $1$ on
$]-\infty,0]$ and vanishing  near $h=+\infty$. We set: 
\begin{equation}\label{triv-infty-F} 
  C = M^{1/2} F''+(1-M)^{1/2}J.
\end{equation}
We are going to check that  $C^2-1\in\cK(\cE'')$. Observe that:
$$ 
   C^2 = M (F'')^2 + 1-M + [M(1-M)]^{1/2}[F'',J] = 1+
    [M(1-M)]^{1/2}[F'',J] \mod \cK(\cE'') 
 $$
where the bracket is  $\ZZ_2$-graded.  
To compute $[F'',J]$, let us proove that $(1+D^2)^{-1/2}$
commutes with $J$. Setting :
 $$
    \Delta = 1+ D^2 = 
             \begin{pmatrix} 
                 \Delta_+ & 0 \\ 
                    0 &  \Delta_-
             \end{pmatrix} 
           =  \begin{pmatrix}
             H_+ & 0 & 0 & 0 \\
             0 & H_- &  0 & 0 \\
             0 & 0 & H_- & 0  \\
             0 & 0 & 0 & H_+
             \end{pmatrix},
 $$
one gets:  
 $$ \Delta J = \begin{pmatrix} 0 & \Delta_+K \\
                        \Delta_- K & 0 \end{pmatrix}, \qquad
 J \Delta = \begin{pmatrix} 0 & K\Delta_- \\
                        K\Delta_+ & 0 \end{pmatrix}, 
 $$
and:
 $$ 
  \Delta_+ K = \begin{pmatrix} 0 & H_+ \\
                        H_- & 0 \end{pmatrix}=K\Delta_-, \qquad
  \Delta_- K = \begin{pmatrix} 0 & H_- \\
                        H_+ & 0 \end{pmatrix}=K\Delta_+, 
 $$
hence $J\Delta=\Delta J$ which implies, using functional calculus, that $J$ commutes with  
$\Delta^{-1/2}$, hence: 
 $$
   [F'',J]=\Delta^{-1/2}(DJ+JD)=\Delta^{-1/2}(2h) =
   (2h)(1+h^2+\lambda^2+ I_2\ot B^2)^{-1/2}.
 $$
Since $h\mapsto M(h)(1-M(h))$ has compact support, we conclude that  
$[M(1-M)]^{1/2}[F'',J]\in \cK(\cE'')$, hence 
$C^2=1 \mod \cK(\cE'')$.

\smallskip We thus get $[\cE'',C]\in K_0(C_0(\RR^2)\ot\cK(L^2(L)))$ and
$C_t = M_t^{1/2} F''+(1-M_t)^{1/2}J$ with
$M_t(h)=M(th)$ provides an operatorial homotopy between $x"=(\cE'',F'')$ and
$(\cE'',C)$. 

\smallskip Using a Fourier transform with 
respect to the variable $\lambda$ and a reparametrization
$\RR\simeq]0,1[$ on $h$, we get an isomorphism 
$\phi: C_0(\RR^2)\ot \cK(L^2(L))\overset{\simeq}{\to} C^*(\cC_L\times TI)$ and $C$ gives
rises to an element still noted $C$ and belonging to
$\Psi^0(\cC_L\times TI,\CC^{4k})$. We now set 
 $$ 
  e=\phi_*(x'')=\phi_*(\cE'',C)=[C^*(\cC_L\times TI,\CC^{4k}),C]\in KK(\cdot, \cC_L\times TI)
 $$
Finally we extend $C$ to $T^qX$ by setting $C= J$ on $TX_+$
thanks to the formula (\ref{triv-infty-F}). Hence:
 $$ 
   i_*(e) = [C^*(T^qX,\CC^{4k}),C]\in KK(\cdot,T^qX) \text{ and } 
   C\in\Psi^0(\Tb,\CC^{4k}),\text{ hence } i_*(e)\in\Delta.
 $$ 

\begin{lem}\label{delta-restricts-onto}
 The equality $\ev[+]([\Delta])= K^{0}(TX_+)$ holds.
\end{lem}
{\bf Proof of the lemma: }
To each $K$-theory class 
$\sigma\in  K^0(T^*X_+)$ we shall associate 
$a_\sigma\in \Delta$ with $(a_\sigma)|_{TX_+}=\sigma$. 
 
\smallskip  Each element of $\sigma\in K^0(T^*X_+)$ can be
represented by a continuous section $f$ over $T^*X_+ \setminus X_+$ of the
bundle $\iso(\pi^*E,\pi^*F)$ for some complex vector bundles $E,F$ over $X_+$ pulled-back by 
$\pi : T^*X_+ \to X_+$. One can assume that $f$ is homogeneous of degre $0$ in the 
fibers of $T^*X_+ $ and independent of $h$ near $\{h= 1\}=\partial X_+$. 

\smallskip One sets $E=F=X_+ \times \CC$ since the general case is identical. Using the
Stone-Weierstrass theorem, we can find $g\in C^\infty(T^*X_+)$ polynomial in 
$\xi$, independent of $h$ near $h=1$ and approximating uniformily $f$ in the corona 
 $\{\xi\in T^*X_+ \ | \ 1/2\le |\xi|\le 2\}$ up to an arbitrary small 
$\varepsilon>0$. Thus, modifying $f$ by:
\begin{equation}\label{analytic-symbol}
 f (x,\xi) = 
 g(x,\frac{\xi}{\sqrt{1+|\xi|^2}})
\end{equation}
one gets another representant of $\sigma\in K^0(T^*X_+)$. 

\smallskip Choosing an exponential map $\theta$ for $\Tb$, a cut-off function $\phi\in
C^\infty_c(\Tb)$ equal to $1$ on units and extending $f$ to $T^*X^{o}$
in the obvious way, one can define the following noncommutative symbol on $X$: 
 \begin{equation}\label{af}  
   a_f(u)(\gamma)\!=\!\mbox{Op}_{\theta,\phi}(f)(u)(\gamma)\! :=\!\int_{ \xi\in
     T^*_{x}X^{o},\  s(\gamma')=x'}\hspace{-1cm}
 e^{i<\theta^{-1}(\gamma'\gamma^{-1}),\xi>}f(x,\xi)
 \phi(\gamma'\gamma^{-1})u(\gamma')d\gamma'd\xi,
\end{equation}    
where $x=r(\gamma)$, $x'=s(\gamma)$ and $u\in C^\infty_c(\Tb)$.
This noncommutative symbol is relatively elliptic on $X$, which is not sufficient
here. We then consider the restriction $a_0$ of $a_f$ at $h=0$. If we note
$\theta_0,\phi_0$ the corresponding restrictions of $\theta,\phi$, we
have, using the same formula as (\ref{af}):
 $$ 
  a_0=\mbox{Op}_{\theta_0,\phi_0}(f_0)\in \cP^0_{inv}(\RR\times L).
 $$
Taking the Fourier transform with respect to the real variable in the
above operator, we get a pseudodifferential operator
$\widehat{a_0}(\lambda)$ on $L$ with parameter $\lambda\in\RR$
which satisfies the condition of ellipticity with parameters, hence,
by a classical result on operators with parameters,
$\widehat{a_0}(\lambda)$ is invertible for large $|\lambda|$. Note
that $\theta,\phi$ can be chosen so that: 
 $$
  \theta_0 : \RR\!\times\! TL\!\ni\!(\lambda,y,V)\!\mapsto\!
  (\lambda,\exp^L_y(V))\!\in \!\RR\!\times\!\cC_L \text{ for small } |V|
  \text{ and }
  \phi_0(\lambda,y,y')=\phi_L(y,y'),
 $$
where $\Exp^L$ is for instance the exponential map associated with the
metric (\ref{product-metric}) and $\phi_L$ is compactly supported in the
range of $\exp_L$ and satisfies $\phi_L(y,y)=1$. It follows that, writing  $x=(h,y)\in [0,1]\times L$; 
$\xi=(\lambda,\eta)\in T_x^*X_+\simeq \RR\times T^*_yL$ and 
$f_0(x,\xi)=f_0(y,\lambda,\eta)$: 
\begin{equation}\label{indicial-family}
   \widehat{a_0}(\lambda)(u)(y)=\int_{y'\in L,\eta\in T^*_yL}
   e^{i<(\exp_y^L)^{-1}(y'),\eta>}f_0(y,\lambda,\eta)\phi_L(y,y')u(y')dy'd\eta
\end{equation}
where $u\in C^\infty(L)$ and $f_0$ denotes the restriction of $f$ at $h=0$. 

\smallskip Observe that $f_0$ has a holomorphic extension 
with respect to the cotangent variable $\lambda\in\RR$ in the strip 
\begin{equation}\label{holomorphic-strip}
   \cB=\{ z=\lambda+iu \in\CC\ | \ -1/2<u <1/2\}.
\end{equation}
Indeed, the following function: 
 $$ 
   f_0(y,z,\eta)=g(x,\frac{(z,\eta)}{\sqrt{1+z^2+|\eta|^2}})
 $$
makes sense as a holomorphic function in $z=\lambda+iu\in\cB$ taking values in
the space $C^\infty(T^*L)$ and is equal to $f_0$ when $u=0$.  
Moreover, for fixed $u\in]-1/2,1/2[$, the function: 
 $$ 
   (y,\lambda,\eta)\mapsto f_0(y,\lambda+iu,\eta)
 $$
is a symbol of order $0$ on $L$ with parameter $\lambda\in\RR$ and one can find a
constant $C$ independent of $y,\lambda,\eta$ such that: 
 $$ 
  |f_0(y,\lambda+iu,\eta)- f_0(y,\lambda,\eta)|\le
  C.u
 $$ 
Since $f_0(y,\lambda,\eta)$ satisfies by construction the condition of
ellipticity for symbols on $L$ with
parameters $\lambda\in\RR$, the previous estimate ensures that the same is true for
$f_0(y,\lambda+iu,\eta)$ assuming that $|u|<\alpha$ for some
$\alpha>0$ small enough. 
In the sequel we restrict the
strip $\cB$ according to this ellipticity condition. 

\smallskip It follows that (\ref{indicial-family}) gives rise to a holomorphic
family $z\mapsto \widehat{a_0}(z)$ taking values in elliptic
pseudodifferential operators on $L$ of order $0$. We have noted
earlier that there exists $z\in\cB$ such that $\widehat{a_0}(z)$ is
invertible, so by a classical result on homolorphic families of
Fredholm operators, the sets:  
 $$ 
   \{ z=\lambda+i\rho\ | \ |\rho|\le\alpha' \text{ and } p_z \text{ is
     not invertible}\}
 $$  
are finite for all $\alpha'<\alpha$. Hence, there exists $\beta$
such that $\widehat{a_0}(\lambda+i\beta)$ is invertible for all
$\lambda\in\RR$. 

\smallskip Observe also that each $\widehat{a_0}(z)$ restricted to horizontal
lines $\mbox{im}(z)=u$ in $\cB$ is a pseudodifferential operator on $L$ with parameter
$\lambda=\mbox{re}(z)$, which allows to define
$a_u\in\cP^*_{inv}(\RR\times L)$ by:
 $$
   \widehat{a_u}(\lambda)=\widehat{a_0}(\lambda+iu)
 $$
Choosing a smooth function
$u(h)$ such that $u(0)=\beta$ and $u(1)=0$, we can  define the
required elliptic noncommutative symbol $a_\sigma$ on $X$ by:
 $$
   a_\sigma|_{X_+}=a_f \text{ and }
   a_\sigma|_h=a_{u(h)}\text{ for all } 0\le h\le 1
 $$  
 
\subsection{Unbounded  noncommutative symbols, Fuchs type noncommutative symbols}\label{FuchsSymbols} 
We can also associate $K$-theory classes to  elliptic noncommutative symbols on
$X$ of positive order. To do that, we state: 
\begin{prop}\label{regular} Let $a\in S^m(X,E)$ be an elliptic
  noncommutative symbol on $X$ with $m>0$.
  Let us consider $a$ as an unbounded operator on $C^*(T^qX,E)$ with
  domain $C^\infty_c(T^qX,E)$. Then its closure $\overline{a}$ is
  regular (\cite{BJ}). 
\end{prop} 
\begin{pf}
The closure of $a$ with
domain $C^\infty_c(\Tb,E)$ is regular as an unbounded operator on
$C^*(\Tb,E)$ and following the proof of this result
in \cite{Vas}, we see that everything 
remains true if we consider $a$ as an unbounded operator on
$C^*(T^qX,E)$ with domain $C^\infty_c(T^qX,E)$. Note that  the
result still holds if $a$ is only relatively elliptic on $X$. 
\end{pf} 
As a consequence, to each elliptic noncommutative symbol $a$ on $X$ of order $m>0$
corresponds a morphism $q(a)=a(1+a^*a)^{-1/2}\in \cL(C^*(T^qX,E))$ and
using the construction of complex powers given in \cite{Vas2}, we get: 
\begin{prop} Let $a\in S^m(X,E)$ be an elliptic noncommutative symbol
 of order $m>0$.\\
1) $(1+a^*a)^{-1/2}$ belongs to $S^{-m}(X,E)$. \\
2) $q(a)$ is an elliptic noncommutative symbol on $X$. 
\end{prop}
\begin{pf} 1) Done in \cite{Vas2}.\\
2) Let $b$ be a parametrix for $a$, that is $ab=1+r$, $ba=1+s$ with
$r,s$ regularizing operators vanishing at $h=0$. Then
$(1+a^*a)^{1/2}b$ is a parametrix for $q(a)$. 
\end{pf}
Now we can associate to each elliptic $a\in S^m(X,E,F)$ the
following $K$-theory class : $[q(a)]$. Note that these
noncommutative symbols do not produce directly {\sl unbounded} $KK$-theoritic
elements (\cite{BJ}) since $(1+a^*a)^{-1/2}$ is not a compact operator on the
$C^*(T^qX)$-Hilbert module $C^*(T^qX,E)$. This defect 
leads us to consider {\sl Fuchs type noncommutative symbols}.  Let $\varphi$ be
a positive smooth increasing function of $h$, equal to $1$ if $h\ge 1$
and satisfying $\varphi(h)=h$ near $h=0$. 
\begin{defn} An element $p\in \Psi^*(T^qX,E)$ is a Fuchs type noncommutative symbol
  on $X$ if $\varphi^lp$ belongs to 
  $S^m(X,E)$ for some $l\in\RR_+$. The infimum of such $l$ is then
  called the fuchs type order of $p$. 
  A Fuchs type noncommutative symbol $p$ with Fuchs type order $l$ strictly positive is
  elliptic if the noncommutative symbol $\varphi^lp$ is elliptic on $X$. 
\end{defn}
For an elliptic Fuchs type noncommutative symbol $p$, we can define as before
$(1+p^*p)^{-1/2}\in \cL(C^*(T^qX,E))$. Thanks to the unbounded
behavior of $p$ with respect to $h$ at $h=0$, the operator
$(1+p^*p)^{-1/2}$ is actually compact so $(C^*(T^qX,E),p)$ provides an
unbounded $KK$-theoritic element in the sense of \cite{BJ}. 

\smallskip Examples of such symbols come from Dirac type operators on $X$, where the latter is provided
with a conical metric $g=dh^2+h^2g_L$, and their typical expression near $h=0$ is: 
 $$ 
   p = h^{-1} \begin{pmatrix} 0 & -\partial_\lambda +S \\ 
   \partial_\lambda +S & 0\end{pmatrix}, \  
 $$ 
where $S$ is a Dirac type operator on $L$. See \cite{DLN2006} for a developpement of this example.

\section{Poincar\'e dual of elliptic noncommutative symbols}
\subsection{Construction of a noncommutative symbol map}
We are going to define a {\sl noncommutative symbol map} for $b$-operators
using a deformation process encoded by $\cG^{q}$. We then get a
generalization of the complete symbol map for manifolds
\cite{Get,Wi1980},  and like the notion that it generalizes, the
noncommutative symbol is not canonical and depends on several
choices: exponential maps, cut-off functions, connections on vector bundles. 
The idea of the construction is very close to \cite{Get,ENN1996}. 
 
We motivate the forthcoming constructions by recalling the case of differential operators on a smooth
manifold $V$. Let $Q$ be a differential operator on $V$ and: 
 $$
  Q(x,D_x)=\sum_\alpha a_\alpha(x)D_x^\alpha
 $$
its expression in a given local chart. For each $t\in ]0,1]$, the differential operator
$P_t$ on $V$ defined locally by: 
\begin{equation}\label{diff-tang-manifold-tpositive}
   P_t(x,D_x)=Q(x,tD_x)
\end{equation}
is well defined and setting: 
\begin{equation}\label{diff-tang-manifold-tzero} 
   P_0(x,D_X)=\sum_\alpha a_\alpha(x)D_X^\alpha\in \mbox{Diff}(T_xV),
\end{equation}
we get a differential operator $P=(P_t)_{t\in [0,1]}$ 
on the tangent groupoid $\cG_V=TV\times\{0\}\cup
\cC_V\times]0,1]$ of $V$. As explained in paragraph
\ref{symbols-manifolds}, $P_0$ represents exactly the (total) symbol
of $Q$. 

\smallskip Let us do the same thing for $b$-differential operators on
$M=\overline{X^{o}}$ with the tangent groupoid $\cG^q$
(\ref{cGq}) of the pseudomanifold $X$. 

\smallskip From now on, the gluing function $l$  (\ref{gluingfunction-l}) is
equal to logarithm function $l=\log$ near $h=0$. Let $Q$ be a $b$-differential operator on
$M$ \cite{Mel}. That means that near $\{h=0\}=\partial M$, one has,
writing $x=(h,y)\in]0,1[\times L$:
 $$ 
  Q = \sum_{k} a_k(h,y,D_y)(h\partial_h)^k 
 $$
where $a_k$ are differential operators on $L$, depending smoothly in $h$.
We define the family $(P_t)_{t\in[0,1]}$ as in
(\ref{diff-tang-manifold-tpositive}), (\ref{diff-tang-manifold-tzero})
on $X_+$, while we set on $X_-$, writing $X=(\lambda,V)\in\RR\times T_yL$:  
\begin{equation}\label{diff-tpositive}
   P_t= \sum_k a_k(h,y,\kappa(t,h)D_y)\left(\frac{t}{l'(h)}\partial_h\right)^k \text{ if }
   t>0, 
\end{equation}
\begin{equation}\label{diff-tzero}  
 P_0= \sum_k a_k(h,y,\tau(h)D_y)(D_\lambda)^k \text{ if }
 h<1 \text{ if } t=0.
\end{equation}
The functions $\kappa$ and $\tau$ are those chosen in
(\ref{gluingfunction-kappa}) and (\ref{gluingfunction-tau}). 
Observe that for $h$ close enough to $0$, (\ref{diff-tpositive}) and (\ref{diff-tzero})
give:
 $$ 
   P_t= \sum_k a_k(h,y,D_y)(th\partial_h)^k \text{ and } P_0= \sum_k a_k(h,y,D_y)(D_\lambda)^k.
 $$
Note that $P_0$ is a noncommutative symbol on $X$ and that $P_0|_{h=0}=p(0,y,D_y,D_\lambda)$ is exactly the {\sl
  indicial operator} of $Q$ \cite{Mel}. Moreover, the full ellipticity of $Q$ as a
$b$-operator (that is its interior ellipticity and the invertibility
of the indicial family) is the same as the ellipticity of $P_0$ as a noncommutative symbol on
$X$. Hence, we have defined a map  
 $$
   \sigma : Q\mapsto P_0
 $$ 
defined on $b$-differential operators and taking values in noncommutative symbols on
$X$. It remains to extend this map to the pseudodifferential case: the
idea is basically the same but things are more technical. 
 
We will use a cover of  $M\times M$ by three open subsets
$R_1,R_2,R_3$ as shown below and a partition of unity
$\omega_1,\omega_2,\omega_3$ subordinated to this cover. 
For instance, $R_1=([0,1/2[\times L)^2$ while 
 $$
  R_2 = \{ (x,x')\in M^2\ | \dist(x,y)<1,\ h(x)+h(x')>3/2 \}
 $$
and $R_3$ is some open neighborhood of the complement of $R_1\cup R_2$ into $M^2$. 

\begin{center} 
\begin{picture}(0,0)%
\includegraphics{partition.pstex}%
\end{picture}%
\setlength{\unitlength}{1782sp}%
\begingroup\makeatletter\ifx\SetFigFont\undefined%
\gdef\SetFigFont#1#2#3#4#5{%
  \reset@font\fontsize{#1}{#2pt}%
  \fontfamily{#3}\fontseries{#4}\fontshape{#5}%
  \selectfont}%
\fi\endgroup%
\begin{picture}(9982,9660)(-959,-8653)
\put(6166,-6316){\makebox(0,0)[lb]{\smash{{\SetFigFont{20}{24.0}{\rmdefault}{\mddefault}{\updefault}{\color[rgb]{0,0,0}$R_3$}%
}}}}
\put(901,-6451){\makebox(0,0)[lb]{\smash{{\SetFigFont{20}{24.0}{\rmdefault}{\mddefault}{\updefault}{\color[rgb]{0,0,0}$R_1$}%
}}}}
\put(5131,-2491){\makebox(0,0)[lb]{\smash{{\SetFigFont{20}{24.0}{\rmdefault}{\mddefault}{\updefault}{\color[rgb]{0,0,0}$R_2$}%
}}}}
\put(811,-871){\makebox(0,0)[lb]{\smash{{\SetFigFont{20}{24.0}{\rmdefault}{\mddefault}{\updefault}{\color[rgb]{0,0,0}$R_3$}%
}}}}
\put(4411,-8521){\makebox(0,0)[lb]{\smash{{\SetFigFont{10}{12.0}{\rmdefault}{\mddefault}{\updefault}{\color[rgb]{0,0,0}$1/2$}%
}}}}
\put(8911,-8386){\makebox(0,0)[lb]{\smash{{\SetFigFont{10}{12.0}{\rmdefault}{\mddefault}{\updefault}{\color[rgb]{0,0,0}$h$}%
}}}}
\put(-944,-3706){\makebox(0,0)[lb]{\smash{{\SetFigFont{10}{12.0}{\rmdefault}{\mddefault}{\updefault}{\color[rgb]{0,0,0}$1/2$}%
}}}}
\put(-314,704){\makebox(0,0)[lb]{\smash{{\SetFigFont{10}{12.0}{\rmdefault}{\mddefault}{\updefault}{\color[rgb]{0,0,0}$h'$}%
}}}}
\end{picture}%

\end{center}
 
\smallskip Let $Q\in\cP^*_b(M)$ with
Schwartz kernel $\kappa$. Let $Q_i$ $i=1,2,3$ be the operators with  Schwartz
kernel  $\kappa_i=\omega_i\kappa$ so that: $Q=Q_1+Q_2+Q_3$. 
 
\smallskip Let us focus on $Q_1$. Applying
a Mellin transfom on $\kappa_1$, we get:
 $$ 
   a_1(h,\eta,y,y') =
   \int_{\RR_+^*}\left(\frac{h}{h'}\right)^{i\eta}
   \kappa_1(h,y,h',y')\frac{dh'}{h'} 
 $$
where $a_1$ is a smooth function of $h\in[0,1/2[$ taking values in the
space of pseudodifferential operators on $L$ with one parameter
$\eta\in\RR$. Using a cut-off function $\phi_1$ such that
$\omega_1\phi_1=\omega_1$, one recovers the action of $Q_1$ on functions
as follows:
 $$
  u\in C^\infty_c(X^{o}),\quad  Q_1 u (h,y) = \int_{\RR_+^*\times\RR} \left(\frac{h}{h'}\right)^{i\eta}
  \left(a_1(h,\eta)\cdot\left(\phi_1(h,h')u(h',.)\right)\right)(y)
  \frac{dh'}{h'}d\eta  
 $$
where $\left(a_1(h,\eta)\cdot u(h',.)\right)(y)$ 
is the action of the operator $a_1(h,\eta)$ on  $u(h',.)\in C^\infty(L)$  evaluated
at $y\in L$.  
Since $\kappa(t,h)=1$ and $l(h)=\log(h)$ when $h\in[0,1/2]$, setting:
 \begin{equation}\label{boundary-term-positive-t}
   u\in C^\infty_c(X^{o}),\quad
   P_{1,t} u (h,y) = \int_{\RR_+^*\times\RR} \left(\frac{h}{h'}\right)^{i\eta/t}
   \left(a_1(h,\eta).\phi_1(h,h')u(h',.)\right)(y)\frac{dh'}{th'}d\eta  
 \end{equation}
for $t>0$ and defining $P_{1,0}\in\Psi^*(\Tb)$ by:
\begin{equation}\label{boundary-term-zero-t}
  u\in C^\infty_c(]0,1/2[\times L\times\RR),\quad
    P_{1,0} u (h,x,\lambda) = \int_{\RR^2} e^{i(\lambda-\lambda').\eta}
           \left(a_1(h,\eta).u(h,.,\lambda')\right)(x)d\lambda'd\eta  
\end{equation}
we get a pseudodifferential operator on $\cG^q$ given by
$P_1=(P_{1,t})_{t\in[0,1]}$ and such that $P_1|_{t=1}=Q_1$.

Since $\kappa_2$ is supported in $R_2$ which is included both in
a compact subset of  $X^{o}\times X^{o}$ and in the range of
$E_{\cG^q}=\Theta^{-1}$, we can set:
 $$ 
   \widehat{a_2}(x,V)=\kappa_2(E_{\cG^q}(x,V,1)) \text{ and } 
    a_2(x,\xi)=\int_{M} e^{i\Theta(x,x',1).\xi}
    \widehat{a_2}(x,\Theta(x,x',1)) d\lambda^1(x').
 $$
Then, choosing any function $\phi_2$ compactly supported in a
neighborhood of $R_{2}$ and satisfying $\omega_2\phi_2=\omega_2$, we have for all functions $u\in C^\infty_c(X^{o})$: 
 $$
   Q_2u(x) = \int_{M\times T^*_xM} e^{i\Theta(x,x',1).\xi}a_2(x,\xi)\phi_2(x,x')u(x')
   d\lambda^1(x')d\xi.
 $$ 
To extend $Q_2$ as we did for $Q_1$, we set:  
\begin{equation}\label{interior-term-positive-t}
   P_{2,t}u(x) =\int_{M\times T^*_xM} e^{i\Theta(x,x',t).\xi}
   a_2(x,\xi)\phi_2(x,x')u(x')d\lambda^t(x')d\xi .
\end{equation}
This is for $t>0$, and we define $P_{2,0}\in\Psi(\Tb)$ by:
\begin{equation}\label{interior-term-zero-t-small-h}
  h<1,\, u\in C^\infty_c(\RR\times L),\, P_{2,0}|_h u(\lambda,y)
 \!\! =\!\!\int_{\RR\times L\times T^*_{(h,y)}M}\hspace{-1.5cm}
  e^{i\Theta(h,\lambda-\lambda',y,y',0).\xi}
   a_2(h,y,\xi)u(\lambda',y')d\lambda^{h,0}(\lambda',y')d\xi     
\end{equation}
and 
\begin{equation}\label{interior-term-zero-t-large-h}
  u\in C^\infty_c(TX_+),\quad P_{2,0}|_{X_+}u(x,X)=
  \int_{T_xM\times T^*_xM}\!\!\!\!\!\!
   e^{i\Theta(x,X-X',0).\xi}a_2(x,\xi) u(x,X')d\lambda^x(X')d\xi.
\end{equation}
The last piece $Q_3$ is smoothing and its Schwartz kernel
$\kappa_3(x,x')$  vanishes both on a neighborhood of the diagonal  and on a
neighborhood of $\partial M\times \partial M$ in $M^2$. This implies
that $\widetilde{\kappa_3}$ defined by  
$\widetilde{\kappa_3}(x,x',t)=\kappa_3(x,x')$
if $t>0$ and $\widetilde{\kappa_3}|_{t=0}=0$, belongs to $C^\infty(\cG^q)$ and 
the behavior of $\kappa_3(x,x')$ near $h(x)=h(x')=0$, resulting from
the assumption that $Q$ is in the small calculus, yields also
$\kappa_3(x,x',t)\in C^*(\cG^q)$. Thus setting for $t>0$ and $u\in C^\infty_c(X^{o})$: 
 $$
  P_{3,t}u(x)=\int_M \kappa_3(x,x')u(x')d\lambda^t(x')  
 $$
and for $t=0$: $P_{3,0}=0$, we have extended $Q_3$ in
$P_3\in\Psi^{-\infty}(\cG^q)\cap C^*(\cG^q)$. We get a linear map: 
 $$ 
   \begin{matrix} 
       \cP_b^*(M) & \longrightarrow & \Psi^*(\cG^q) \\ 
          Q  & \longmapsto &  P_Q:=P_1+P_2+P_3 .
   \end{matrix}
 $$
Restricting $P_Q$ at $t=0$ gives the desired noncommutative symbol map:
\begin{defn} Let $Q\in\cP^*_b(M)$. With
  the notations above, we define the noncommutative symbol of $Q$ by : 
  $$ 
    \sigma(Q) = P|_{t=0}=P_{1,0}+P_{2,0} \in S^*(X). 
  $$
\end{defn} 
The following facts are obvious: 
\begin{rem} 
  \begin{enumerate} 
    \item  $Q$ is fully elliptic as a $b$-operator if and only if
      $\sigma(Q)$ is elliptic as a noncommutative symbol on $X$. 
    \item If $P\in\cP^p_b(M)$ and $Q\in\cP^q_b(M)$ then 
      $$ 
         \sigma(PQ)=\sigma(P)\sigma(Q) \text{ modulo } S^{p+q-1}(X).
      $$
    \item Everything above can be written in the same way for
      operators acting on sections of a vector bundle. 
  \end{enumerate}
\end{rem}
The noncommutative symbol map depends on the choices of the cover of $M^2$, of the
partition of unity and of the exponential map of $\cG^q$, but the
$K$-theory class of the noncommutative symbols of fully elliptic $b$-operators does
not depend on these choices, as we will see in the next paragraph. 

\smallskip Conversely, one can define a {\sl quantification} map $\op_b$ which is a quasi
inverse of $\sigma$. We describe it now. Let us choose $\omega\in C^\infty_c([0,1/2[)$ such that
$\omega(h)=1$ near $h=0$. Let $a$ be a noncommutative symbol on $X$ and write 
 $$ 
   a=\omega a + (1-\omega)a = a_1+a_2.
 $$
We extend $a_1$ as a $\cG^q$-pseudodifferential operator
$\widetilde{a_1}$ by reverting the process used in
(\ref{boundary-term-positive-t},\ref{boundary-term-zero-t}). 
Let $f_2$ be a symbol of $a_2$ viewed as a
pseudodifferential operator on $T^qX$. That means that $f_2$ is an
ordinary symbol on $A^*(T^qX)=T^*X^{o}$ such that: 
  $$
    a_2 = \op_{T^qX}(f_2) \text{ modulo } S^{-\infty}_0(X,E) 
  $$
where $\op_{T^qX}$ is given by:
 $$ 
   u\in C^\infty_c(T^qX), \ \op_{T^qX}(f_2)(u)(\gamma)=
   \int_{(T^qX)_{s(\gamma)}\times T^*_{r(\gamma)}X^{o}}
  \hspace{-1.5cm}e^{i<\Theta^q(\gamma'\gamma^{-1}),\xi>}
  f_2(r(\gamma),\xi)\phi(\gamma'\gamma^{-1})u(\gamma')d\lambda^{s(\gamma)}d\xi .
 $$
Here $\Theta^q=(E_{T^qX})^{-1}$ is the inverse of the
exponential map of $T^qX$ given by restriction of $E_\cH$, and $\phi$ is a cut-off function equal to $1$ on units
and supported in the range of $\Theta^q$. 
We can use the formulae (\ref{interior-term-positive-t},\ref{interior-term-zero-t-small-h},
\ref{interior-term-zero-t-large-h}) to build from $f_2$ a   
 $\cG^q$-pseudodifferential operator $\widetilde{a_2}$ with
the property :
 $$
    \widetilde{a_2}|_{t=0} = a_2 \text{ modulo } S^{-\infty}_0(X). 
 $$
Thus we get an approximate lifting of noncommutative symbols: 
\begin{equation}\label{liftsymbolstocGq}
\widetilde{a}:=\widetilde{a_1}+\widetilde{a_2}\in\Psi^*_b(\cG_q)
\end{equation} 
satisfying:  
 $$
   \widetilde{a}|_{t=0}=a\text{ modulo } S^{-\infty}_0(X),
 $$
\begin{equation}\label{opb}
  \op_b(a):= \widetilde{a}|_{t=1}\in\cP^*_b(M).
\end{equation}
By construction: 
\begin{equation}\label{rightinverse-symbol}
   \sigma(\op_b(a))=a  \text{ modulo } S^{-\infty}_0(X)
\end{equation}
In the same way, if $P\in\cP_b^*(M)$ then $\op_b(\sigma(P))-P$ is
a smoothing operator with vanishing indicial operator.  
\subsection{The Poincar\'e duality as a noncommutative symbol map}
All ingredients are now at hands to finish. Observe that to each
fully elliptic $b$-operator  $P: C^\infty(M,E)\to C^\infty(M,F)$
acting on sections of complex vectors bundles $E,F$, corresponds a $K$-homology class 
$[P]=[(L^2(M;E\oplus F;d\lambda^1),\rho_1,\mathbf{P})]$ on $X$ where: 
\begin{itemize}
  \item $L^2(M;E\oplus F;d\lambda^1)$ is the $\ZZ_2$-graded Hilbert space modeled on the measure
    $d\lambda^1$  and on  product type hermitian structures on $E$ and
    $F$,
  \item $\rho_1$ is the action of $C(X)$ onto $L^2(M;E\oplus F;d\lambda^1)$ in the
    natural way through the quotient map $M\to X=M/\overline{X_-}$, 
  \item $\mathbf{P}=\begin{pmatrix} 0& Q\\ P & 0 \end{pmatrix}$ where
    $Q$ is a full parametrix of $P$.
\end{itemize}
\begin{thm} With the previous notations and those of the definition
  \ref{bounded-symbol}, the isomorphism $\Sigma^q:K_0(X)\to
  K^0(T^qX)$ defined in (\ref{PD-map-conic-q-case}) is given
  by: 
\begin{equation}\label{interpretationPD}
    [P] \longmapsto [\sigma(P)] 
\end{equation}
\end{thm}
\begin{rem}
 \item[$\bullet$] Recall that from theorem \ref{symbgenerate}, we know that
   every $K$-theory class $[a]\in K^0(T^qX)$ has a representant $a$ among
   elliptic noncommutative symbols on $X$. From (\ref{rightinverse-symbol}),  we
   know that $a$ is in the same $K$-theory class than the noncommutative symbol
   $\sigma(P)$ of a fully elliptic $b$-operator $P$. Eventually,
   since  $\Sigma^q$ is an isomorphism,  we get that each
   $K$-homology class of $X$ is represented by a fully
  elliptic $b$-operator.  

 \item[$\bullet$] Using the deformation of $T^{q}X$ into $T^{c}X$, leading to
   the $KK$-equivalence $T^qX\sim T^cX$, one could also get a
   concrete interpretation for $\Sigma^c$.  However, the adequate adaptation
   of the notion of noncommutative symbols is more difficult to relate
   directly to what is done in boundary values problems or former
   studies about pseudodifferential calculus for groupoids. 
\end{rem}
\begin{pf} Let $P:C^\infty(M,E_0)\to C^\infty(M,E_1)$ be a fully
  elliptic $0$-order b-operator. Let 
$E=E_0\oplus E_1$ and $a=\sigma(P)\in S^0(X,E)$. We need to prove that 
$\Sigma^q[P]=[a]$ or equivalently that
$[P]=[a]\underset{T^qX}{\otimes} D^q$ (cf. section \ref{section2}). 
Recall that: 
\begin{equation}\label{recallPDformula}
   [a]\underset{T^qX}{\otimes} D^q=\mathbf{s}_X([a])\ot\Phi^q \ot \partial^q.
\end{equation}
Firstly, $\mathbf{s}_X([a])\ot\Phi^q \in KK(X,T^qX)$ is represented by: 
\begin{equation}\label{firststep}
    (C^*(T^qX,E), \rho, \mathbf{a})
\end{equation}
where $\rho : C(X)\To \cL(C^*(T^qX,E))$ is given by 
$\rho(f)(\xi)(\gamma)=\xi(\gamma)f(\pi^q(\gamma))$.

\smallskip The next step is to find 
\begin{equation}\label{secondstep}
   [\widetilde{\cE},\widetilde{\rho},\widetilde{\mathbf{a}}]\in KK(X,\cG^q)
\end{equation}
such that 
\begin{equation}\label{pb-lift-sigma}
   (e^q_0)_*[\widetilde{\cE},\widetilde{\rho},\widetilde{\mathbf{a}}]=[C^*(T^qX,E), \rho, \mathbf{a}].
\end{equation}
The desired lifting is made as follows. Let us note again $E$ the pull
back of the original bundle $E$ to $X^{o}\times [0,1]$ with the range 
map of $\cG^q$. Let $\widetilde{\pi^q}$  be the composite map of the
range map $\cG^q\to X^{o}\times[0,1]$  with the projection maps 
$X^{o}\times[0,1]\to X^{o}$ and $X^{o}\to X=X^{o}/\overline{X_-}$. 

\smallskip We set $\widetilde{\cE}=C^*(\cG^q,E)$, we define $\widetilde{\rho}$ by
$\widetilde{\rho}(f)(\xi)(\gamma)=f(\widetilde{\pi^q}(\gamma))\xi(\gamma)$,
and $ \widetilde{\mathbf{a}}$ is defined from $\mathbf{a}$ using
(\ref{liftsymbolstocGq}). By construction
$\widetilde{\ba}\in\Psi^0_b(\cG^q,E)\subset \cL(C^*(\cG^q,E))$ and 
$\widetilde{\ba}^2=1$ modulo 
$\Psi^{-1}_{b,0}(\cG^q,E)\subset \cK(C^*(\cG^q,E))$.
It follows that the triple 
$(\widetilde{\cE},\widetilde{\rho},\widetilde{\mathbf{a}})$ defined above
satisfies (\ref{secondstep}) and (\ref{pb-lift-sigma}). 
Evaluating this element at $t=1$ gives: 
 $$
  (e^q_1)_*[\widetilde{\cE},\widetilde{\rho},\widetilde{a}]
    =[\widetilde{\cE}|_{t=1},\widetilde{\rho}_{t=1},\op_b(a)] \in
      KK(X,\cC_{X^{o}})
 $$
and applying the Morita equivalence $C^*(\cC_{X^{o}})\overset{\nu}{\sim} \CC$
produces the final result: 
 $$
   [a]\underset{T^qX}{\otimes} D^q= [\widetilde{\cE}|_{t=1},\widetilde{\rho}_{t=1},\op_b(a)]
                \otimes \nu =
                [L^2(M;E;d\lambda^1),\rho_1,P]=[P]\in
                KK(X,\cdot)=K_0(X)
 $$
\end{pf}


\subsection{Index map}
Since $X$ is a compact Hausdorff space the
map $p$ sending $X$ to a point gives rise to a morphism: 
 $$
   p_*\ : K_0(X) \To K_0(\cdot)=\ZZ
 $$
called, for obvious reasons, the index map. We can capture $p_*$ with
the pre-Dirac element and Poincar\'e duality:
\begin{prop}\label{analytic-index}
Let us denote by $\hbox{Ind}^q$ the map:
 $$ 
  \hbox{Ind}^q : [a]\in K(T^qX) \mapsto [a]\otimes \partial^q\in\ZZ,
 $$
then the following holds: 
 $$ 
   \forall [a]\in K(T^qX),\qquad p_*((\Sigma^q)^{-1}[a]) =\hbox{Ind}^q[a] .
 $$
In other words, the index of a fully elliptic $0$-order $b$-operator 
$P: C^\infty(M,E)\to C^\infty(M,F)$ viewed as a Fredholm operator
between $L^2(M,E,d\lambda^1)$ and $L^2(M,F,d\lambda^1)$ is equal to 
$\Ind[\sigma(P)]{q}$.
\end{prop}
\begin{pf} The homomorphism $\CC\to C(X)$ corresponding to
  $p:X\to\cdot$ is denoted by $\widetilde{p}$. 
We have   
  $$ p_*((\Sigma^q)^{-1}[a]) = [\widetilde{p}]\ot \left([a]\underset{T^qX}{\otimes}([\Phi^q]\otimes
  \partial^q)\right)= [\widetilde{p}]\ot\mathbf{s}_{X}[a]\ot [\Phi^q]\ot\partial^q. 
  $$
  Observe that 
 $$
 [\widetilde{p}]\ot\mathbf{s}_{X}[a]=[\widetilde{p}]\underset{\CC}{\ot}[a]$$
so by the commutativity of the Kasparov product over $\CC$:
 $$[\widetilde{p}]\ot\mathbf{s}_{X}[a]= \ot[a]\underset{\CC}{\ot}
 [\widetilde{p}]= a\ot \mathbf{s}_{T^{q}X}[\widetilde{p}]. $$
But $\mathbf{s}_{T^{q}X}[\widetilde{p}]\ot[\Phi^q]$ is equal to the class of the
identity homomorphism of $C^{*}(T^{q}X)$, hence:
\begin{eqnarray*}
  p_*((\Sigma^q)^{-1}[a])&=& ([\widetilde{p}]\ot\mathbf{s}_{X}[a])\ot
  [\Phi^q]\ot\partial^q \\
 &=& (a\ot \mathbf{s}_{T^{q}X}[\widetilde{p}])\ot
 [\Phi^q]\ot\partial^q \\
 &=& a\ot(\mathbf{s}_{T^{q}X}[\widetilde{p}]\ot
 [\Phi^q])\ot\partial^q \\
&=& a\ot\partial^q 
\end{eqnarray*}

\end{pf}

\bibliographystyle{plain} 
\bibliography{biblio.bib} 
 
\end{document}